\documentclass[a4paper]{article}

\usepackage{amsfonts}
\usepackage{amsopn}
\usepackage{amsmath}
\usepackage{amsthm}
\usepackage{latexsym}
\usepackage{amssymb}
\usepackage{graphicx}

\newcommand{\p}{\partial}

\renewcommand{\phi}{\varphi}

\newcommand{\comp}{^{\sf c}}
\newcommand{\R}{\mathbb R}

\newcommand{\cA}{{\mathcal A}}

\newcommand{\cC}{{\mathcal C}}

\newcommand{\cK}{{\mathcal K}}

\newcommand{\cN}{{\mathcal N}}

\newcommand{\rE}{{\mathrm E}}

\newcommand{\dd}{\hspace{1pt}{\mathrm d}}

\newcommand{\ee}{{\mathrm e}\hspace{1pt}}

\begin{document}

\title{Random set solutions to elliptic and hyperbolic partial differential equations}

\author{Jelena Karaka\v{s}evi\'{c}\thanks{Unit of Engineering Mathematics, University of Innsbruck,
Technikerstra\ss e 13, 6020 Innsbruck,
Austria, (jelena.karakasevic.mail@gmail.com)}
\and
Michael Oberguggenberger\thanks{Unit of Engineering Mathematics, University of Innsbruck,
Technikerstra\ss e 13, 6020 Innsbruck,
Austria, (michael.oberguggenberger@uibk.ac.at)}
}

\date{}
\maketitle

\begin{abstract}
The past decades have seen increasing interest in modelling uncertainty by heterogeneous methods, combining probability and interval analysis, especially for assessing parameter uncertainty in engineering models. A unifying mathematical framework admitting the combination of a wide range of such methods is the theory of random sets, describing input and output of a structural model by set-valued random variables. The purpose of this paper is to highlight the mathematics behind this approach. The modelling and computational implications are discussed and demonstrated with the help of prototypical partial differential equations---a scalar elliptic equation from elastostatics and hyperbolic systems arising in elastodynamics.
\end{abstract}


\section{Introduction}
\label{sec:intro}
Many models in structural engineering are cast in the form of (systems of) ordinary or partial differential equations. Static problems, e.g. from elasticity theory, give rise to elliptic partial differential equations. Dynamic problems (e.g. vibrations of elastic bodies or acoustic wave propagation) typically are formulated in terms of hyperbolic partial differential equations or systems. The quantities of interest are solutions to (elliptic) equilibrium equations or to (hyperbolic) equations of motion. These equations contain coefficients that---as a rule---depend on further parameters, such as material constants or geometric dimensions. In addition, the parameters may vary in space and time. An important task in engineering design is the analysis of the uncertainty of the computed output. This uncertainty traces back to uncertainty in the input parameters, which may be due to many sources: random fluctuations, lack of information, conflicting assessments, random measurement errors, systematic measurement errors, fluctuations due to spatial inhomogeneity, errors made by assigning parameter status to state variables, model insufficiency, just to name a few. Uncertainty in the model itself is important, but will be left aside in this paper.

Available information on parameter variability may consist in frequency distributions obtained from large samples, values from small samples or single measurements, interval bounds, experts' point estimates, or educated guesses from experience.

Accordingly, an increasing desire has risen in the engineering community to introduce methods of quantifying the uncertainty beyond probability, such as
interval analysis, set-valued models, fuzzy sets, evidence theory, random sets, sets of probability measures, imprecise probability, lower and upper previsions, info-gap-analysis, etc. This paper focusses on combinations of probability and interval analysis.
This includes random variables, random fields, intervals, interval fields, random variables with interval parameters and random fields with interval parameters (such as mean, variance, correlation length).

As already noted in \cite{2006:alvarez,2015:oberguggenberger,2008:MOFellin}, all these types and their combinations can be accommodated in the framework of \emph{random sets}. Random sets are set-valued random variables, whereby the values may be sets of real numbers, subsets of finite dimensional vector spaces, or subsets of function spaces. For example, random variables are single-valued random sets and intervals are random sets on a one-point probability space. The purpose of this paper is to highlight some important mathematical modelling questions that arise when using random sets.
The considerations will be exemplified with the help of a prototypical elliptic equation and a class of hyperbolic systems.

When uncertainties in the coefficient functions of ordinary or partial differential equations are modelled as random sets, random fields with interval parameters, or (finite) interval fields, the solutions of these equations are no longer single-valued, but rather random sets with values in function spaces. Evaluating the solution at a fixed point in space or time will result in a real-valued random set. However, here an important issue arises: in order to properly qualify as a random set, the set-valued solution has to satisfy a certain measurability condition (to be detailed in Section \ref{subsec:rs}). It turns out that the essential ingredient for this condition to hold is that the solution should depend continuously on the coefficients of the equation, and in turn also on the hyperparameters of the random variables or random fields defining the coefficients. This continuous dependence is already needed in pure interval modelling to prove that the values of the solutions are intervals. It is also needed in probabilistic modelling to infer that the solutions are random variables or random fields. In these two (classical) cases, the arguments are rather straightforward. However, in models combining interval and probabilistic inputs, additional subtleties are involved. One of the main goals of the present article is to formulate a rigorous mathematical framework that allows one to infer that the solutions are random sets. The elliptic and the hyperbolic case require different arguments, the essence of which will be worked out in the corresponding sections.

Here is a short review of relevant literature on uncertainty quantification in engineering models, combining interval methods with probabilistic methods. Adopting the language of \cite{2018:jiang}, the structural model can be cast in the form $Z = g(X,Y)$ where $X$ and $Y$ refer to different types of inputs in a hybrid probability-interval uncertainty analysis; $Z$ is the system response or output, $g$ is the performance function or input-output map. The map $g$ can be given, for example, by an ordinary or partial differential equation, a finite element model, or a black box algorithm. Here $X$ refers to the probabilistic part, $Y$ to the interval part---both may be multi-dimensional. If both parts are present, the output is necessarily both random and set-valued. There are three subcases.  In the first case, $X$ and $Y$ are considered as separate, non-interacting entities. The set-valued output arises through their combination in the input-output map $g$. In the second case, $X = X(Y)$, that is, $X$ is a random variable or random field with interval hyperparameters. This is often referred to as an imprecise random variable or an imprecise random field. The third possibility is that $Y = Y(X)$, the interval bounds are random variables. This is usually referred to as a random interval, a special case of a random set.

Typical representative articles concerning the first case are \cite{2016:muscolino,2022:Sofi}, where $n$-DOF structures with interval material parameters under stochastic excitation have been studied. Structures described by stochastic ordinary differential equations with random interval parameters have been addressed in \cite{2010:schmelzer:adam}, while wave equations with interval propagation speed under space-time white noise excitation have been addressed in \cite{2019:wurzer}. Combinations of intervals and random vectors as input parameters in beams and truss structures have been considered in \cite{2010:gao}. The paper \cite{2018:lue} merges intervals and random variables by means of evidence theory (finite random sets) in an application to squeal analysis of a brake system.

The second case extends from imprecise random variables as material parameters in frame and grid structures \cite{2020:Sofi} to finite element models with imprecise random field inputs \cite{2016:do,2019:faes}, employing a Karhunen-Lo\`eve expansion with interval coefficients. The applications in Sections 4 and 5 of this paper are in the same spirit. A combination of inputs consisting of imprecise random variables and imprecise random fields can be found in \cite{2015:oberguggenberger}.

As mentioned, the third case concerns random interval input and is addressed in the finite dimensional version in \cite{2010:gao,2008:MOFellin}, in combination with stochastic differential equations in the quoted paper \cite{2010:schmelzer:adam} and in \cite{2010:schmelzer,2013:schmelzer}.

There are various papers which contrast and compare interval modelling and probabilistic modelling side by side in beam and plate structures \cite{2013:muscolino,2015:sofi,2020:Sofi,2018:sofi:romeo}. The question of identification and model calibration under hybrid uncertainty has been addressed in the literature as well. A procedure for identifying interval models from measurements has been proposed in \cite{2017:faes}. Constructing random sets from statistical data has been considered in \cite{2008:MOFellin}. In the recent paper \cite{2022:kitahara} a Bayesian-type updating procedure for random variables with interval constraints has been presented, while a whole range of calibration methods has been combined in \cite{2022:gray} for quantification of mixed uncertainties, reliability analysis and design of structures.

Survey papers on modelling uncertainty in engineering by hybrid probability-interval approaches are \cite{2013:beer:ferson:kreinovich,2013:oberguggenberger,2018:jiang}. A recent review paper on numerical methods for uncertainty propagation in probability-box models is \cite{2021:FaesDaub}. Further details and references on numerical algorithms can be found in Subsection \ref{subsec:num}. For the theory of random sets, their interpretations and applications the reader is referred to the monographs \cite{1977:castaing,2005:molchanov,2006:nguyen}, notably \cite{2010:bernardini:book} for engineering applications.

The plan of the paper is as follows. Section 2 will be devoted to the required mathematical set-up. These facts are known, but it appears useful to elaborate them in some detail for the later applications. Section 3 singles out the principal mathematical argument that allows one to prove that the output quantity is a random set. The combination of arguments giving a full proof of this result seems to be new. Section 4 details the considerations for elliptic equations as an application of the findings of Section 3, including a simple explicit numerical example. Section 5 deals with hyperbolic systems and elaborates the continuity argument in the hyperbolic case, also having a short look at hyperbolic equations and systems in higher space dimensions. Section 6 contains a comparison with the parametric approach and numerical aspects, as well as some remarks on modelling correlations. Detailed proofs can be found in the dissertation of the first author \cite{2020:nedeljkovic}.

The present paper is an extended and updated version of the proceedings contribution \cite{2021:MOKarakasevic}.

\section{Mathematical preliminaries}
\label{sec:math}

This section specifies some mathematical notions which are required for a rigorous treatment of random sets. Recall that a subset $S$ of a finite dimensional vector space such as $n$-dimensional Euclidean space $\R^d$ is \emph{open}, if for each point $x$ in $S$, each ball of the form $\{y\in\R^d:\|y-x\|<\eta\}$ with center in $x$ and sufficiently small radius $\eta$ is also contained in $S$. The set $S$ is \emph{closed}, if with each converging sequence whose terms belong to $S$, the limit also belongs to $S$. The subset $S$ is \emph{bounded}, if the norms $\|x\|$ of the elements of $S$ are all smaller than a fixed constant. Further, finite dimensional Euclidean spaces are \emph{separable}, that is, they contain a countable dense subset (a sequence of points such that every ball in $\R^d$ contains at least one element of the sequence), and they are \emph{complete} (that is, if $x_k$ is a sequence such that the sum of its norms $\sum_{k=1}^\infty\|x_k\|$ converges, then the sum $\sum_{k=1}^\infty x_k$ also converges).

All these notions can be literally generalized to infinite dimensional normed spaces. A typical example is the space of continuous functions on a bounded, closed interval $I$, denoted by $\cC(I)$, with norm $\|f\| = \max_{x\in I}|f(x)|$. This space is separable and complete.

\subsection{Random variables}
\label{subsec:rv}

A real-valued random variable $a$ is characterized by its distribution function $F_a(x)$ which defines the probabilities $P$ of basic events of the form $\{a \le x\}$ via
$P(a \leq x) = F_a(x)$. This point of view suffices for most purposes. However, in the sequel \emph{families of random variables} will be needed, and it will be necessary to specify the domains of the random variables under consideration. The usual mathematical setting is a probability space $(\Omega,\Sigma, P)$ where $\Omega$ is a set (the collection of all elementary events), $\Sigma$ is the family of measurable subsets of $\Omega$, and $P$ is a probability measure (that assigns to every set belonging to $\Sigma$ a number between 0 and 1, satisfying the usual axioms). In this point of view, a random variable is a mapping from $\Omega$ into the real line such that the events $\{\omega\in\Omega: a(\omega) \le x\}$ are measurable (i.e., belong to $\Sigma$). Thus a probability $P\big(\{\omega\in\Omega: a(\omega) \le x\}\big)$ can be assigned to such events, denoted again by $P(a \le x)$, and these probabilities in turn define the distribution function of the random variable $a$. From the probabilities $P(a \le x)$, the probabilities of open and closed intervals, their complements, their countable unions and intersections etc. can be computed. The family of subsets of the real line which can be obtained by indefinite repetition of these set theoretic operations defines the so-called \emph{Borel $\sigma$-algebra}.
The Borel $\sigma$-algebra can also be characterized as the smallest $\sigma$-algebra containing all open subsets of $\R$.
A random variable has the property that all sets of the form $\{\omega\in\Omega: a(\omega) \in B\}$, where $B$ is any Borel set, are measurable. Accordingly, random variables are often referred to as measurable functions.

If the probability space is finite, say $\Omega = \{1,\ldots, m\}$, one may take the power set of $\Omega$ as $\sigma$-algebra, and every real-valued function on $\Omega$ is measurable. A typical infinite probability space is standard Gaussian space with $\Omega = \R$, $\Sigma$ the Borel measurable subsets of $\R$, and the probability of an event $B$ given by the Gaussian integral
$
P(B) = (2\pi)^{-1/2}\int_B\ee^{-\omega^2/2}\dd\omega;
$
the probability measure $P$ is given through the Gaussian density $\ee^{-\omega^2/2}/\sqrt{2\pi}$. On standard Gaussian space, every continuous real-valued function on $\Omega$ is a random variable.

These notions can be easily extended to finite dimensional random variables, i.e., with values in $\R^d$ (then the joint probability distributions have to be specified).  Measurable functions with values in an infinite dimensional normed space can be defined in the same way.

\subsection{Random fields}
\label{subsec:rf}

A random field is a process that assigns a random variable $q(x)$ to every point $x$ in a region $D$ in space. To define the probabilistic properties of the field, the joint distributions of the values at any finite number of points $q(x_1), \ldots q(x_n)$ should be specified.
If the random field is Gaussian, it is completely specified by the mean value $\mu_q = \rE(q(x))$ and the second moments, i.e., the covariance $C(x,y) ={\mathrm{COV}}(q(x),q(y))$ for any two points $x,y$. If it is (weakly)  homogeneous and isotropic, the covariance depends only on the distance $\rho = \|x-y\|$ of the points and is of the form
$
    C(x,y) = \sigma^2 c(\rho)
$
with the field variance $\sigma^2$ and the autocorrelation function $c(\rho)$. A typical autocorrelation function is of the form
\begin{equation}
\label{eq:ACF}
    c(\rho) = \exp\Big(-\frac{|\rho|}{\ell}\Big),
\end{equation}
where $\ell$ is the so-called correlation length. The standard method of simulation of a random field is based on the \emph{Karhunen-Lo\`eve expansion}. If $q(x)$ is any mean zero random field with finite second moments given by the covariance $C(x,y)$, its Karhunen-Lo\`eve expansion is obtained by solving the eigenvalue problem
\begin{equation}
\label{eq:eigenvalueProblem}
   \int_D C(x,y)\phi_k(y)\dd y = c_k\phi_k(x)
\end{equation}
which has a sequence of positive eigenvalues $c_k$ and orthonormal eigenfunctions $\phi_k(x)$ (orthonormality in mean square).
Then
\begin{equation}
\label{eq:KL}
   q(x) = \sum_{k=1}^\infty\sqrt{c_k}\xi_k\phi_k(x)
\end{equation}
where the $\xi_k$ are uncorrelated random variables with unit variance. If the process is Gaussian, the $\xi_k$ are independent and distributed according to $\cN(0,1)$.

For the numerical simulation, the spatial region is discretized by a grid and the $\phi_k$ are taken, e.g., piecewise constant on the grid elements. The eigenvalue problem becomes a matrix eigenvalue problem, and the series (\ref{eq:KL}), with approximate eigenvalues and eigenfunctions, truncated after a finite number $M$ of terms,
can be used for Monte Carlo simulations of the field trajectories (i.e., realizations of the field).

In case $q(x)$ is a mean zero homogeneous Gaussian random field on $D = \R$ with autocorrelation function (\ref{eq:ACF}) and variance $\sigma^2$, an alternative to the Karhunen-Lo\`eve expansion is the representation as an Ornstein-Uhlenbeck process, namely as solution to the Langevin stochastic differential equation
\begin{equation}\label{eq:Langevin}
\dd q = - \frac1{\ell} q\,\dd x  + \sqrt{\frac{2}{\ell}}\,\sigma\,\dd W, \quad q|_{x=0}\sim\cN(0,\sigma^2),
\end{equation}
see e.g. \cite{1974:arnold}. Here $W = W(x,\omega)$ denotes Wiener process on the real line.
In this case, realizations of the random field can be conveniently generated by simulating solutions to the Langevin equation. An additional advantage of the representation through (\ref{eq:Langevin}) is that the dependence of $q(x)$ on the field parameters $\sigma, \ell$ is more explicit than through the Karhunen-Lo\`eve expansion (where $\sigma, \ell$ enter in the coefficients $c_k$ and the eigenfunctions $\varphi_k$). The representation (\ref{eq:Langevin}) allows one to directly prove that the assignment $(\sigma,\ell)\to q(x,\omega)$ is continuous at fixed $x$ and $\omega$, an important ingredient in the random set case discussed in Section~\ref{sec:parrf}.

\subsection{Random sets}
\label{subsec:rs}

In general, a {\em random set} is a set-valued random variable satisfying certain measurability conditions, to be detailed below.
The simplest case arises when the underlying probability space is finite. In this case,
one speaks of {\em finite random sets} or {\em Dempster-Shafer structures}.
Such a structure is given by finitely many closed subsets
$A_i, i = 1, \dots, m $ of a target space $\cA$, usually Euclidean space $\R^d$, called the {\em focal
elements}, each of which comes with a {\em probability weight} $p_i = p(A_i)$, $\sum p_i = 1$.

Viewed as a random set, a Dempster-Shafer structure is given by an $m$-point probability space
$\Omega = \{1, 2, \ldots, m\}$ with probability masses $\{p_1, p_2,\ldots, p_m\}$; the assignment $i\to A_i$ is the defining
set-valued random variable $A$.
However, starting with Subsection~\ref{subsec:rsf} and particularly in Section~\ref{sec:parrf}, random sets defined on infinite probability spaces will be needed. Also, the target space can be infinite-dimensional, for example, a function space containing the system outputs.

In the sequel, $(\Omega,\Sigma,P)$ denotes a general probability space. The target space will be either Euclidean space $\R^d$ or more generally a normed, complete and separable space $\cA$. In analogy to the Borel $\sigma$-algebra on $\R$ introduced in Subsection~\ref{subsec:rv}, the Borel $\sigma$-algebra on $\cA$ is defined as the smallest $\sigma$-algebra containing all open subsets of $\cA$.
Given a set-valued map
$\omega\to A(\omega)$ where each $A(\omega)$ is a subset of the target space $\cA$, upper and lower inverses of subsets $B$ of the target space $\cA$ are defined by
\begin{equation}\label{inverses}
   A^-(B) = \{\omega\in\Omega: A(\omega) \cap B \neq \emptyset\},\qquad A_-(B) = \{\omega\in\Omega: A(\omega) \subset B \}.
\end{equation}
The requirements that $\omega\to A(\omega)$ is a \emph{random set} are
\begin{itemize}
\item each $A(\omega)$ is closed;
\item the upper inverses $A^-(B)$ are measurable subsets of $\Omega$ for every closed subset $B$ of the target space $\cA$.
\end{itemize}
In the later applications, one would like to admit also open sets $B$ or more generally Borel sets $B$ in the second condition. In addition, it will be useful to relate it to the existence of certain measurable selections.
Recall that a \emph{measurable selection} of a random set $A$ is a random variable $a$ such that $a(\omega) \in A(\omega)$ for all $\omega \in \Omega$. A \emph{Castaing representation} is a sequence of measurable selections $a_k(\omega) \in A(\omega)$ such that for all $\omega$, $\{a_k(\omega):k = 1,2,3,\ldots\}$ is dense in $A(\omega)$.

In this respect, an important technical tool is the \emph{fundamental measurability theorem}. In order to obtain all subsequent equivalences, one needs to assume that $(\Omega,\Sigma,P)$ is a \emph{complete probability space}, meaning that every subset of a set of probability zero belongs to $\Sigma$. This can always be achieved by enlarging $\Sigma$ so as to include these subsets.

The fundamental measurability theorem states that, assuming the first condition (each $A(\omega)$ closed) that the following conditions are equivalent:
\begin{itemize}
\item the upper inverses $A^-(B)$ are measurable subsets of $\Omega$ for every closed subset $B$ of the target space $\cA$;
\item the upper inverses $A^-(B)$ are measurable subsets of $\Omega$ for every open subset $B$ of the target space $\cA$;
\item the upper inverses $A^-(B)$ are measurable subsets of $\Omega$ for every Borel subset $B$ of the target space $\cA$;
\item $\omega\to A(\omega)$ admits a Castaing representation.
\end{itemize}
For the rather elaborate proof, the reader is referred to \cite{1977:castaing,2005:molchanov}.
Here are some first consequences of the fundamental measurability theorem.
First, the lower inverses of events $B$ can be written as $A_-(B) = \{\omega\in\Omega: A(\omega) \cap B\comp \neq \emptyset\}\comp$, where $B\comp$ denotes the complement of $B$, and hence are measurable due to the fundamental measurability theorem and the measurability of $A^-(B)$.
Second, thanks to the required measurability properties, it is legitimate to introduce upper and lower probabilities
\begin{eqnarray}
\overline{P}(B) &=& P\big(A^-(B))\ = \ P(\{\omega\in\Omega: A(\omega) \cap B \neq \emptyset\})\label{eq:uppprob}\\
\underline{P}(B) &=& P\big(A_-(B))\ = \ P(\{\omega\in\Omega: A(\omega) \subset B \})\label{eq:lowprob}
\end{eqnarray}
for any Borel set $B\subset \cA$.

The concept of random sets encompasses many other well-established methods of uncertainty modelling. Clearly, every interval defines a random set (as a random variable on a one-point probability space). Similarly, every random variable $a$ on a probability space $(\Omega,\Sigma,P)$ defines a (single-valued) random set $\omega\to A(\omega) = \{a(\omega)\}$. The upper and lower probability of an event $B$ coincide with its probability $P(B)$. Further,
every normalized fuzzy set can be viewed as a random set. Here the probability space is the interval $\Omega = [0,1]$, equipped with the uniform probability distribution, and the focal elements $A(\omega)$ are just the $\omega$-level sets. It is not difficult
to prove that the possibility measure $\pi(B)$ of a subset $B$ of the real line coincides with its upper probability \cite{2002:goodman}. Indeed, and rather obviously,
\[
  \pi(B) = \sup_{\omega\in[0,1]}\{A(\omega)\cap B \neq \emptyset\} = P(\{\omega\in\Omega: A(\omega) \cap B \neq \emptyset\}) = \overline{P}(B).
\]
In the special case of finite random sets, the measurability condition holds automatically; the only remaining condition is that the focal elements $A_i$ are closed.

\subsection{Random sets generated by families of random variables}
\label{subsec:rsf}

\begin{figure}[htb]
\begin{center}
\includegraphics[width = 7cm]{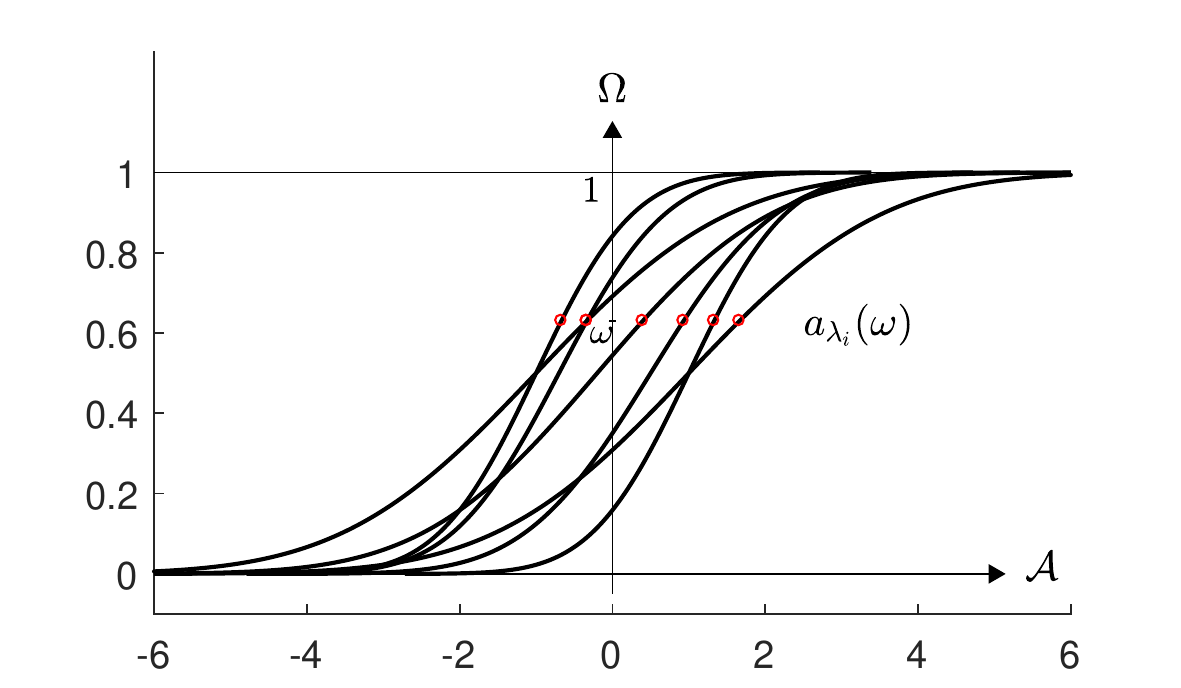}\qquad\qquad \includegraphics[width = 7cm]{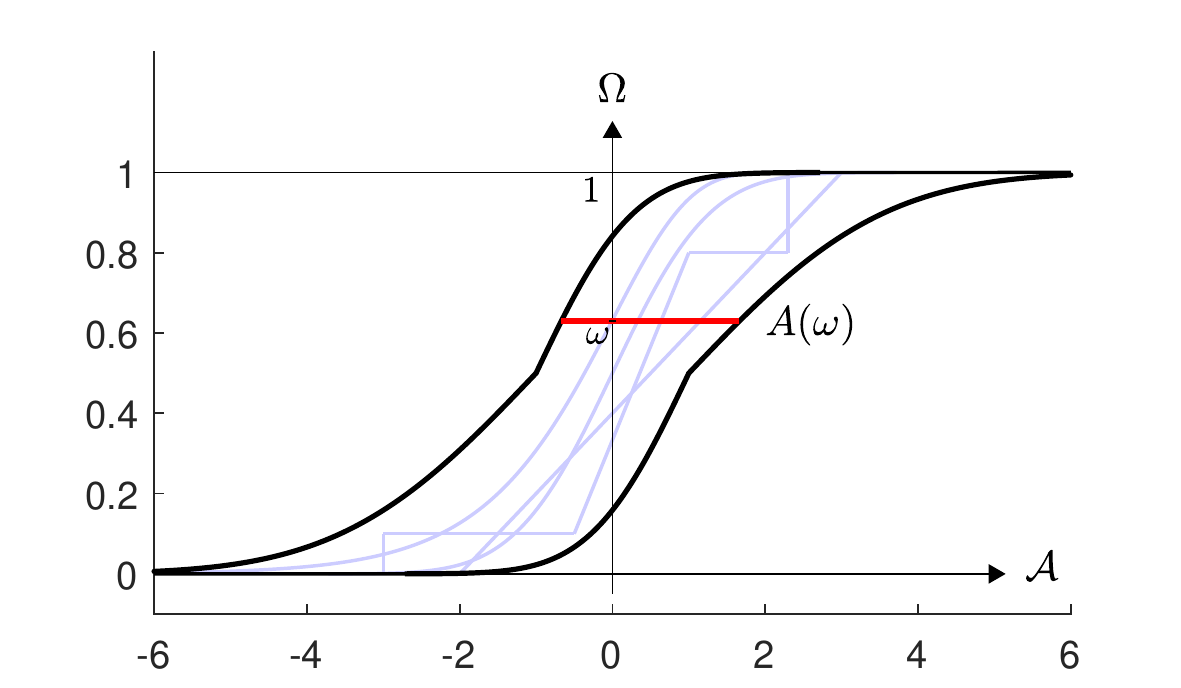}
\caption{Some members of the family $a_\lambda$ (left), random set with focal element (right). The parameters $\lambda = (\mu,\sigma)$ were taken as $-1\leq\mu\leq 1$, $1\leq\sigma\leq 2$.}
\end{center}
\label{fig:1}
\end{figure}

To set the stage, consider random variables $a_\lambda$, defined on the same probability space $(\Omega,\Sigma,P)$, and depending on a set of parameters $\lambda\in\Lambda$.
Suppose the random variables have values in a target space $\cA$. Then one can define a set-valued map
$
    \omega \to A(\omega) =  \{a_\lambda(\omega):\lambda\in\Lambda\}.
$
The question is whether this assignment defines a random set. Before answering this question, consider the example of an imprecise Gaussian family
\begin{equation}\label{eq:GaussRS}
   a_\lambda(\omega) = \mu + \sigma \Phi^{-1}(\omega)
\end{equation}
with interval parameters
$
\lambda =(\mu,\sigma)\in
   \Lambda = [\underline{\mu},\overline{\mu}]\times [\underline{\sigma},\overline{\sigma}].
$
Here the underlying probability space is the unit interval $\Omega = (0,1)$ with the uniform distribution, and $\Phi$ is the standard normal distribution function.
It follows that each random variable $a_\lambda$ is distributed according to ${\mathcal N}(\mu,\sigma^2)$. The corresponding subset of the real line is
\[
   A(\omega) = [\underline{a}(\omega),\overline{a}(\omega)]
\] with
\[
 \underline{a}(\omega) = \min\big(a_\lambda(\omega):\lambda\in\Lambda\big),\quad
      \overline{a}(\omega) = \max\big(a_\lambda(\omega):\lambda\in\Lambda\big).
\]
Due to the fact that the parameter set $\Lambda$ is a (two-dimensional) closed and bounded interval and the fact that the assignment $(\mu,\sigma)\to \mu + \sigma \Phi^{-1}(\omega)$ is continuous at fixed $\omega$, it follows that
each $A(\omega)$ is a closed and bounded interval, indeed. Some members of the family and the resulting random intervals are depicted in Figure \ref{fig:1} (left).

Anticipating that $\omega\to A(\omega)$ is a random set, one may define the upper and lower distribution functions
\begin{equation}\label{eq:lucdfs}
\underline{F}(b) = P(\omega: \overline{a}(\omega)\leq b), \quad \overline{F}(b) = P(\omega: \underline{a}(\omega)\leq b).
\end{equation}
Anticipating also that the functions $\underline{a}$, $\overline{a}$ bounding the random intervals are random variables, one obtains simply
that $\underline{F}$ is the distribution function of $\overline{a}$ and $\overline{F}$ is the distribution function of $\underline{a}$.

Here is the argument why $\omega\to A(\omega)$ is a random set. First, one constructs a Castaing representation by taking a countable dense subset of parameter values
$\lambda_k = (\mu_k,\sigma_k)$ in
$\Lambda = [\underline{\mu},\overline{\mu}]\times [\underline{\sigma},\overline{\sigma}]$, for example, all points with rational coordinates.
Again from the continuity of the assignment $\lambda\to a_\lambda(\omega)$ at fixed $\omega$, i.e., $(\mu,\sigma)\to\mu + \sigma \Phi^{-1}(\omega)$, see (\ref{eq:GaussRS}), it follows that the set $\{a_{\lambda_k}(\omega): k = 1,2,3,\ldots\}$ is dense in $A(\omega)$ for every fixed $\omega$. Let $B$ be an open set. Then
\begin{equation} \label{eq:meas}
\begin{array}{rcl}
   A^-(B) &=& \{\omega\in\Omega: A(\omega) \cap B \neq \emptyset\}\\
        &=& \{\omega\in\Omega: \mbox{ there\ is\ } k \mbox{\ such\ that\ } a_{\lambda_k}(\omega)\in B\}\\
        &=&  \bigcup_k \{\omega\in\Omega: a_{\lambda_k}(\omega)\in B\}
\end{array}
\end{equation}
is measurable as a countable union of measurable sets; note that the individual sets are measurable, because the maps $\omega\to a_\lambda(\omega)$ are measurable at each fixed $\lambda$. The fundamental measurability theorem implies that $A^-(B)$ is measurable for every Borel set $B$, in particular, for every closed set $B$. This is the required property for $\omega \to A(\omega)$ to be a random set.

Using again that $\{a_{\lambda_k}(\omega): k = 1,2,3,\ldots\}$ is dense in $A(\omega)$, it follows that $\underline{a}(\omega) = \inf\{a_{\lambda_k}(\omega): k = 1,2,3,\ldots\}$ and $\overline{a}(\omega) = \sup\{a_{\lambda_k}(\omega): k = 1,2,3,\ldots\}$
are measurable as infima and suprema of countably many random variables.

\emph{Comparison with the concept of $p$-boxes.} Given two cumulative distribution functions $\underline{F} \leq \overline{F}$ on the real line, a \emph{probability box} or \emph{$p$-box}, as introduced in \cite{2003:ferson}, is a collection of cumulative distribution functions $F$ such that $\underline{F}(b) \leq F(b) \leq \overline{F}(b)$ for all $b\in\R$. The random set $\omega\to A(\omega)$ above can be viewed as a $p$-box through (\ref{eq:lucdfs}), but it collects only the cumulative distribution functions corresponding to the random variables $a_\lambda$, $\lambda\in\Lambda$. Such $p$-boxes are often referred to as \emph{distributional} or \emph{parametric $p$-boxes}. The argument above shows that distributional $p$-boxes are random sets, provided the underlying family of random variables depends continuously on the parameter $\lambda$. On the other hand, $p$-boxes that are understood as containing all cumulative distribution functions between $\underline{F}(b)$ and $\overline{F}(b)$ (the original concept of \cite{2003:ferson}) may be referred to as \emph{distribution-free} or \emph{non-parametric $p$-boxes}. One may define a corresponding set-valued map by
\[
  \omega\to [\overline{F}^{(-1)}(\omega), \underline{F}^{(-1)}(\omega)]
\]
employing the quasi-inverses of $\underline{F}^{(-1)}$ and $\overline{F}^{(-1)}$. It has been shown in \cite{2006:alvarez} that this assignment defines a random set as well; see also the discussion in \cite{2018:alvarez}. The distributional and the distribution-free $p$-boxes corresponding to $\underline{F}(b)$ and $\overline{F}(b)$ in (\ref{eq:lucdfs}) look the same, but the distribution-free version contains many probability distribution functions other than the one coming from the Gaussian family (indicated in Figure \ref{fig:1}, right).

\section{Random sets generated by parametrized random fields}
\label{sec:parrf}

This central section lays down the basic structure that will allow one to prove that set-valued maps obtained as solutions to partial differential equations whose coefficients are given as intervals, random fields, or parametrized random fields, are indeed random sets.

\subsection{The general measurability result}
\label{subsec:meas}
For the present purpose, the most general situation to be considered is the following.
\begin{itemize}
\item $(\Omega,\Sigma,P)$ is a complete probability space;
\item the target space $\cA$ is a normed, complete and separable space;
\item $\Lambda$ is a bounded and closed subset of some Euclidean space;
\item $a_\lambda(\omega)$, $\lambda\in \Lambda$ is a family of random variables with values in $\cA$, that is, the maps $\omega\to a_\lambda(\omega)$ are measurable for each $\lambda$;
\item the random variables depend continuously on $\lambda$, that is, the maps
$
\Lambda \to \cA: \lambda\to a_\lambda(\omega)
$
are continuous for each $\omega$.
\end{itemize}
Then the following assertions hold: The set-valued map
\[
    \omega \to A(\omega) =  \{a_\lambda(\omega):\lambda\in\Lambda\}
\]
defines a random set. All focal elements $A(\omega)$ are bounded and closed subsets of $\cA$. In addition, if $\Lambda$ is a (multi-dimensional) interval and $\cA = \R$, the sets $A(\omega)$ are intervals.

\emph{Proof of assertion.}
First, it is clear from the continuity assumption that all $A(\omega)$ are bounded and closed subsets of $\cA$. The assertion that they are intervals under the additional specializing hypothesis follows from the continuity as well. To prove the random set property, one has to show that the upper inverses $A^-(B)$ are measurable for every
open subset $B$ of $\cA$. This follows exactly by the argument given in equation (\ref{eq:meas}) and applying the fundamental measurability theorem.

\emph{Example.}
Consider a function $u = u(\lambda,a)$ that depends on an interval parameter $\lambda$, say $\lambda \in [\underline{\lambda},\overline{\lambda}]$, and a fixed random variable $a$. Such a situation arises, for example, when $u$ is the response of a system with mixed interval and random uncertainty. Both $\lambda$ and $a$ can be multi-dimensional, so this covers the case of discretized interval and random fields. Assume that the dependence of $u$ on $\lambda$ and $a$ is continuous.
The functions $a_\lambda(\omega) = u(\lambda, a(\omega))$ form a family of random variables exactly of the type discussed here. Collecting the responses in
$
   A(\omega) = \{u(\lambda, a(\omega)):\lambda\in\Lambda\}
$
therefore leads to a random set, that is actually a random interval $\omega\to[\underline{a}(\omega),\overline{a}(\omega)]$ if the response is one-dimensional.

\subsection{Continuous dependence of random fields on their parameters}
\label{subsec:contrf}

A more involved instance of a random set is obtained by taking a random field on the real line ($x\in\R$) whose parameters are intervals. For example, consider mean zero random fields $q(x,\omega,\ell)$  with fixed variance $\sigma^2$, but correlation length $\ell$ varying in an interval $[\underline{\ell},\overline{\ell}]$. A common probability space for these random fields with different correlation lengths is needed. A convenient way of generating such random fields, in the case of the autocorrelation function (\ref{eq:ACF}), is by means of the Langevin equation (\ref{eq:Langevin}). The common probability space is Wiener space $(\Omega_W,\Sigma_W,P_W)$ (see Appendix).

Due to the fact that the random fields $q(x,\omega,\ell)$ are solutions to (\ref{eq:Langevin}), the following properties hold; see e.g. \cite{2010:schmelzer,2013:schmelzer}:
\begin{itemize}
\item For fixed $x\in\R$ and $\ell\in[\underline{\ell},\overline{\ell}]$, the map $\omega \to q(x,\omega,\ell)$ is measurable.
\item For fixed $x\in\R$ and $\omega\in\Omega_W$, the map $\ell \to q(x,\omega,\ell)$ is continuous.
\end{itemize}
Thus at every chosen point $x$, the assignment
$\omega \to Q(x,\omega) = \{q(x,\omega,\ell): \ell \in [\underline{\ell},\overline{\ell}]\}$
defines a (scalar) random set, as can be seen by the same arguments as in the example above. The trajectories of the random field are interval-valued curves, given by $x\to [\min(q(x,\omega,\ell),\max(q(x,\omega,\ell)]$, where the minimum and maximum are obtained by letting $\ell$ run through $[\underline{\ell},\overline{\ell}]$.

Alternatively, the Karhunen-Lo\`eve expansion (\ref{eq:KL}) can be used to place the random fields $q(x,\omega,\ell)$ in a common probability space $\Omega$, not depending on $\ell$. Formally, this can be done by considering an infinite product of standard Gaussian spaces with random elements $\omega = (\omega_1,\omega_2,\omega_3,\ldots)$. The required sequence of Gaussian variables $\xi_k$ is then given by $\xi_k(\omega) = \omega_k$ (see Appendix), and the Karhunen-Lo\`eve expansion reads
\begin{equation}
\label{eq:KLpar}
   q(x,\omega,\ell) = \sum_{k=1}^\infty\sqrt{c_{k,\ell}}\xi_k(\omega)\phi_{k,\ell}(x).
\end{equation}
In general, there is no way to assert that the eigenvalues $c_{k,\ell}$ and eigenfunctions $\phi_{k,\ell}(x)$ depend continuously on the correlation length $\ell$, even after discretization and reduction to a matrix eigenvalue problem. Thus one has to analyze the situation in each case individually. The case of the autocovariance function (\ref{eq:ACF}) on a finite interval can be treated explicitly, thanks to the results in \cite{1991:ghanem}. Without restriction of generality, one may assume that the finite interval is $D = [-1,1]$. Then the eigenvalue problem (\ref{eq:eigenvalueProblem}) becomes
\begin{equation}
\int_{-1}^1{\ee^{-|x-y|/\ell}\phi(y)}\,\dd y=c\phi(x).
\label{eq:KLint}
\end{equation}
Following \cite{1991:ghanem}, problem (\ref{eq:KLint}) can be solved in the following way. First, the two equations
\begin{equation}
\label{eq:tan}
\frac{1}{\ell}-\alpha\tan{\alpha}=0,\quad{\rm respectively}\quad \alpha+\frac{1}{\ell}\tan{\alpha}=0
\end{equation}
have two sequences of positive solutions denoted by $\alpha_k = \alpha_k(\ell)$ and $\alpha^\ast_k = \alpha^\ast_k(\ell)$, for even and odd $k$, respectively. The resulting eigenvalues and eigenfunctions are
\begin{equation}
\label{eq:eigenv}
c_{k,\ell}=\frac{2\ell}{1+\ell^2\alpha_k^2(\ell)}, \quad c_{k,\ell}^\ast=\frac{2\ell}{1+\ell^2\alpha^{\ast 2}_{k}(\ell)}
\end{equation}
\begin{equation}
\label{eq:eigenf}
\phi_{k,\ell}(x)=\frac{\cos{\alpha_k(\ell) x}}{\sqrt{1+\frac{\sin{2\alpha_k(\ell)}}{2\alpha_k(\ell)}}},\quad
\phi^\ast_{k,\ell}(x)=\frac{\sin{\alpha^\ast_k(\ell) x}}{\sqrt{1-\frac{\sin{2\alpha^\ast_k(\ell)}}{2\alpha^\ast_k(\ell)}}}
\end{equation}
for even and odd $k$, respectively.
Therefore, the one-dimensional random field $q(x,\omega,\ell)$ with covariance function given by the equation (\ref{eq:ACF}) has the Karhunen-Lo\`eve expansion
\begin{equation}
\label{eq:KLprocess}
q(x,\omega,\ell)=\lim_{M\to\infty}\sum\limits_{k=1}^M\Big(\sqrt{c_{k,\ell}}\xi_k(\omega)\phi_{k,\ell}(x) + \sqrt{c_{k,\ell}^\ast}\xi^\ast_k(\omega)\phi^\ast_{k,\ell}(x)\Big).
\end{equation}
For the truncated Karhunen-Lo\`eve expansion, one may take a $2M$-dimensional standard Gaussian space as probability space $\Omega$, on which the random variables $\xi_k$ and $\xi^\ast_k$ are defined.
It is clear that the solutions $\alpha_k$ and $\alpha_k^\ast$ of the deterministic equations (\ref{eq:tan}) depend continuously on $\ell$. The same is true for the eigenvalues as well as eigenfunctions shown in (\ref{eq:eigenv}), (\ref{eq:eigenf}).
Thus truncating the Karhunen-Lo\`eve expansion (\ref{eq:KLprocess}) at a finite $M$, the resulting random field depends continuously on the correlation length $\ell$ at fixed $\omega$.

\subsection{Random fields as random variables valued in function spaces}
\label{subsec:functrf}

Consider the random field $q(x,\omega,\ell)$ as in Subsection \ref{subsec:contrf}. The \emph{paths or trajectories (in space)} of the random field are the maps $x\to q(x,\omega,\ell)$, with fixed $\omega$ and $\ell$. It is well-known from the theory of stochastic differential equations that the paths are continuous \cite{1974:arnold}. This also follows from the Kolmogorov-Chentsov theorem \cite{1997:kallenberg}, because the autocovariance function (\ref{eq:ACF}) satisfies the hypotheses of this theorem. However, for the purposes of the present paper, more is needed.

First, jointly continuous dependence on $x$ and $\ell$ is required. This is easy to achieve, because one can consider the \emph{paths in space and parameter set} as the maps $(x,\ell)\to q(x,\omega,\ell)$, with fixed $\omega$ and apply the previous arguments.
Indeed, when using a truncated Karhunen-Lo\`eve expansion, the argument was just given after (\ref{eq:KLprocess}). When the Langevin equation is used, one may resort to the results of \cite{2010:schmelzer,2013:schmelzer} which assert the continuous dependence, jointly in $x$ and $\ell$.

Second, the random fields have to be interpreted as random variables in the space of continuous functions. This is because they enter as variable coefficients in the partial differential equations to be solved, and the solution depends on the whole trajectories, not only the values of the random fields at individual points $x$. Let $x$ belong to an interval $I$ and $\ell$ to an interval $L$, both closed and bounded. Continuity of the paths in $x$ and $\ell$ can be stated by saying that the map $q(\omega):(x,\ell)\to q(x,\omega,\ell)$ belongs to the space of continuous function $\cC(I\times L)$ for fixed $\omega$. Thus the random field induces a map $\Omega \to \cC(I\times L)$, $\omega\to q(\omega)$. Standard arguments from the theory of stochastic processes show that this map is measurable, that is, a random variable with values in $\cC(I\times L)$. For details of the proof, see \cite[Section 1.2]{2020:nedeljkovic}.
The argument also shows that
\begin{itemize}
\item for fixed $\ell\in L$, the map $\Omega \to \cC(I)$, $\omega \to q(\omega,\ell)$ is measurable;
\item for fixed $\omega\in\Omega$, the map $L \to \cC(I)$, $\ell \to q(\omega,\ell)$ is continuous.
\end{itemize}
Note that the same random field $q$ may be considered in three different ways as a measurable function of $\omega$: as an element of the function space $\cC(I\times L)$, denoted by $q(\omega)$, as an element of the function space $\cC(I)$ at fixed correlation length $\ell$, denoted by $q(\omega,\ell)$, and as a real number at each fixed $x$ and $\ell$, denoted by $q(x,\omega,\ell)$. Whenever its values are collected in a random set, the corresponding capital letter $Q$ will be used.

\subsection{The mean values of a parametrized random field as an interval field}
\label{subsec:means}
Let $\omega\to A(\omega)$ be a random set whose focal elements $A(\omega) = [\underline{a}(\omega),\overline{a}(\omega)]$ are closed and bounded intervals. Let $\{a_k(\omega)$, k = 1,2,3,\ldots\} be a Castaing representation.
The \emph{Aumann expectation} $\rE(A)$ is defined as the closure of the set of all expectation values of the measurable selections $a_k$ for $k = 1,2,3,\ldots$, that is, $\rE(A) = {\rm cl}\{\rE(a_k): k = 1,2,3,\ldots\}$, see \cite{2005:molchanov}. It is clear that in the considered case, the Aumann expectation is simply the interval
$
   \rE(A) = [\rE(\underline{a}), \rE(\overline{a})].
$
Next, consider a random set generated by a parametrized random field $\omega \to Q(x,\omega) = \{q(x,\omega,\ell): \ell \in [\underline{\ell},\overline{\ell}]\}$ as in Subsection \ref{subsec:contrf}. Thanks to the continuous dependence on the parameter $\ell$ discussed in Subsection \ref{subsec:contrf} each $Q(x,\omega)$ is a closed and bounded interval at fixed $x$. One may form the Aumann expectation $\rE(Q(x))$, which is a deterministic interval at fixed $x$. Thus the assignment $x\to \rE(Q(x))$ defines an interval field.

\section{Application to elliptic equations}
\label{sec:elastostatics}

A prototypical stochastic elliptic partial differential equation, arising e.g. in elastostatics, is the boundary value problem
\begin{equation}\label{eq:ellprob}
\begin{array}{rclrcl}
	-\,{\rm div\,}\big(a(x)\,{\rm grad\,}u(x)\big)&=&f(x),&x&{\rm in}&D\\
	u(x)&=&0,&x&{\rm in}&\partial D\\
\end{array}
\end{equation}
as discussed, e.g., in \cite{2010:grigoriu}.
Here $D$ is a bounded open domain in $\R^d$ and $\p D$ is its boundary. For example, in two space dimensions $u$ might be the transversal displacement of a non-uniform membrane under transversal load $f$. The spatially varying elastic properties are subsumed in the coefficient $a(x)$. Alternatively, in three space dimensions, $u(x)$ could be the pressure of a fluid in a porous medium with permeability $a(x)$ and source rate $f(x)$ \cite{2008:matthies}.

This section serves to demonstrate how the random field methods can be applied in practice. The required solution theory is collected, and a simple numerical example is elaborated.

In the foundations of the finite element method, problem (\ref{eq:ellprob}) is commonly solved by variational methods. Suitable function spaces are defined as follows: $L^2(D)$ is the space of square integrable functions; $H^1(D)$ is the space of square integrable functions all whose first order partial derivatives are square integrable as well; $H^1_0(D)$ is the subspace of functions that vanish on the boundary $\p D$. Both $L^2(D)$ and $H^1_0(D)$ are normed spaces; the squares of the norms are given by
\[
    \|f\|_{L^2(D)}^2 = \int_D |f(x)|^2\dd x, \quad \|u\|_{H^1(D)}^2 = \|u\|_{L^2(D)}^2 + \|{\rm grad\,u}\|_{L^2(D)}^2.
\]
The spaces are complete and separable. The variational formulation of equation (\ref{eq:ellprob}) then is:

\begin{quote}
Find $u\in H^1_0(D)$ such that
\begin{equation}\label{eq:varprob}
   \int_D a(x)\,{\rm grad\,} u(x)\cdot{\rm grad\,} v(x)\dd x + \int_D f(x)v(x)\dd x = 0
\end{equation}
for all $v\in H^1_0(D)$.
\end{quote}

\emph{The deterministic case.}
Assume that $a$ is a continuous function on $D$ which is bounded from above and below, that is, $a\in \cC(D)$ and $\alpha \leq a(x) \leq \beta$ for some constants $\alpha, \beta > 0$ and all $x$ in $D$. Further, assume that $f$ is square integrable on $D$. Under these conditions, it is well-known that the variational problem (\ref{eq:varprob}) has a unique solution $u$ in $H^1_0(D)$, \cite{1967:necas,1975:treves}. For the purpose of this example, the source term $f$ will be fixed. The solution $u$ depends on the coefficient function $a$; this dependence will be made explicit in the notation $u = u_a$. Denote the set of functions which are continuous on $D$ up to the boundary with values between $\alpha$ and $\beta$ by $\cC(D; \alpha, \beta)$. It is equipped with the norm $\|a\| = \max_{x\in D}|a(x)|$. The crucial assertion is that the map
\[
   \cC(D; \alpha, \beta)\to H^1_0(D): a \to u_a
\]
is continuous. This non-trivial result can be proved by invoking the variational formulation (\ref{eq:varprob}), as it was done in \cite{2020:nedeljkovic}. A different proof, for more general partial differential equations, can be found in \cite{1967:necas}.

\emph{The random field case.}
To start with, let $a(x,\omega)$, $\omega\in\Omega$, be a random field with continuous paths lying between $\alpha$ and $\beta$. As noted in Subsection \ref{subsec:functrf}, it can be viewed as a random variable with values in $\cC(D; \alpha, \beta)$. Due to the continuity of the above map, the assignment $\omega\to u_{a(\omega)}$ is measurable, thus the solution $u_a$ is a random variable with values in $H^1_0(D)$. Intuitively, this means that $u_a$ is a random field whose paths belong to $H^1_0(D)$. In order to interpret it as a classical random field, one should be able to assign values $u_{a(\omega)}(x)$ at each point $x$ in $D$. In one space dimension, this is easy because $H^1_0(D)$ is continuously imbedded in $\cC(D)$, and so the pointwise evaluations
\[
   \cC(D; \alpha, \beta)\to \R: a \to u_a(x)
\]
are continuous as well. Thus the maps $\omega\to u_{a(\omega)}(x)$ are random variables for each $x$, and the paths $x\to u_{a(\omega)}(x)$ are continuous functions. In space dimensions two and three, more theory is needed. First, the input random field $a(x,\omega)$ should have Lipschitz continuous paths, that is, satisfying an inequality of the form $|a(x,\omega) - a(y,\omega)|\leq C|x-y|$ for some constant $C$ and all $x,y$. This is not true of the random field generated by the Langevin equation, but it is true of the random field generated by the truncated Karhunen-Lo\`eve expansion. Then one can invoke regularity theory saying that the solution to problem (\ref{eq:ellprob}) is in $H^2(D')$ in any open subregion $D'$ of $D$ and depends continuously on $a$ in this space \cite{1967:necas}. Again, $H^2(D')$ is continuously imbedded in $\cC(D')$, and so one obtains a classical random field $u_{a(\omega)}(x)$, defined for all $x$ in the interior of $D$.

\emph{The random set case.}
Let $a(x)$ be a parametrized random field. More specifically,
it will be constructed as follows. In two space dimensions, take a closed and bounded interval $I$ such that $D\subset I\times I$ and a random field of the form $q(x,\omega,\ell)$ as in Subsection \ref{subsec:contrf} and set
\begin{equation}\label{eq:coeff}
   a(x,\omega,\ell) = \mu(x) + \chi\big(q(x_1,\omega,\ell)q(x_2,\omega,\ell)\big), \quad {\rm for} \quad x = (x_1,x_2){\ \rm in\ } D.
\end{equation}
Here $\mu(x) > 0$ is the mean field (assumed fixed and deterministic) and $\chi$ is a cut-off function; this is needed because the individual realizations of a Gaussian random field can be arbitrarily large with non-zero probability. In the end, the cut-off has to be chosen so as to guarantee that $\alpha \leq a(x,\omega,\ell) \leq \beta$.

One could also use different correlation lengths in the two variables and/or multiply the fields with different standard deviations $\sigma_1,\sigma_2$. The construction in one or three space dimensions is analogous.

When $\ell$ varies in an interval $L$, a set-valued solution to (\ref{eq:ellprob}) can be defined as $U(\omega) = \{u_{a(\omega,\ell)}:\ell\in L\}$. Combining the continuity of the map $a\to u_a $ with the statement at the end of Subsection \ref{subsec:functrf}, one obtains that
\begin{itemize}
\item for fixed $\ell\in L$, the map $\Omega \to H^1_0(D)$, $\omega \to u_{a(\omega,\ell)}$ is measurable;
\item for fixed $\omega\in\Omega$, the map $L \to H^1_0(D)$, $\ell \to u_{a(\omega,\ell)}$ is continuous.
\end{itemize}
It follows from Subsection \ref{subsec:meas} that $U(\omega)$ is a random set in the target space $H^1_0(D)$. Again, in one space dimension one may form the set-valued functions $U(\omega,x) = \{u_{a(\omega,\ell)}(x):\ell\in L\}$ at each $x$ in $D$. By the same argument, the $U(\omega,x)$ are random sets in $\R$; they form the point evaluations of the random set solution $U(\omega)$ at the points $x$ in $D$. In space dimensions two and three, a modified argument as in the random field case has to be applied.

\begin{figure}[htb]
  \centering
  \includegraphics[width = 5.16cm]{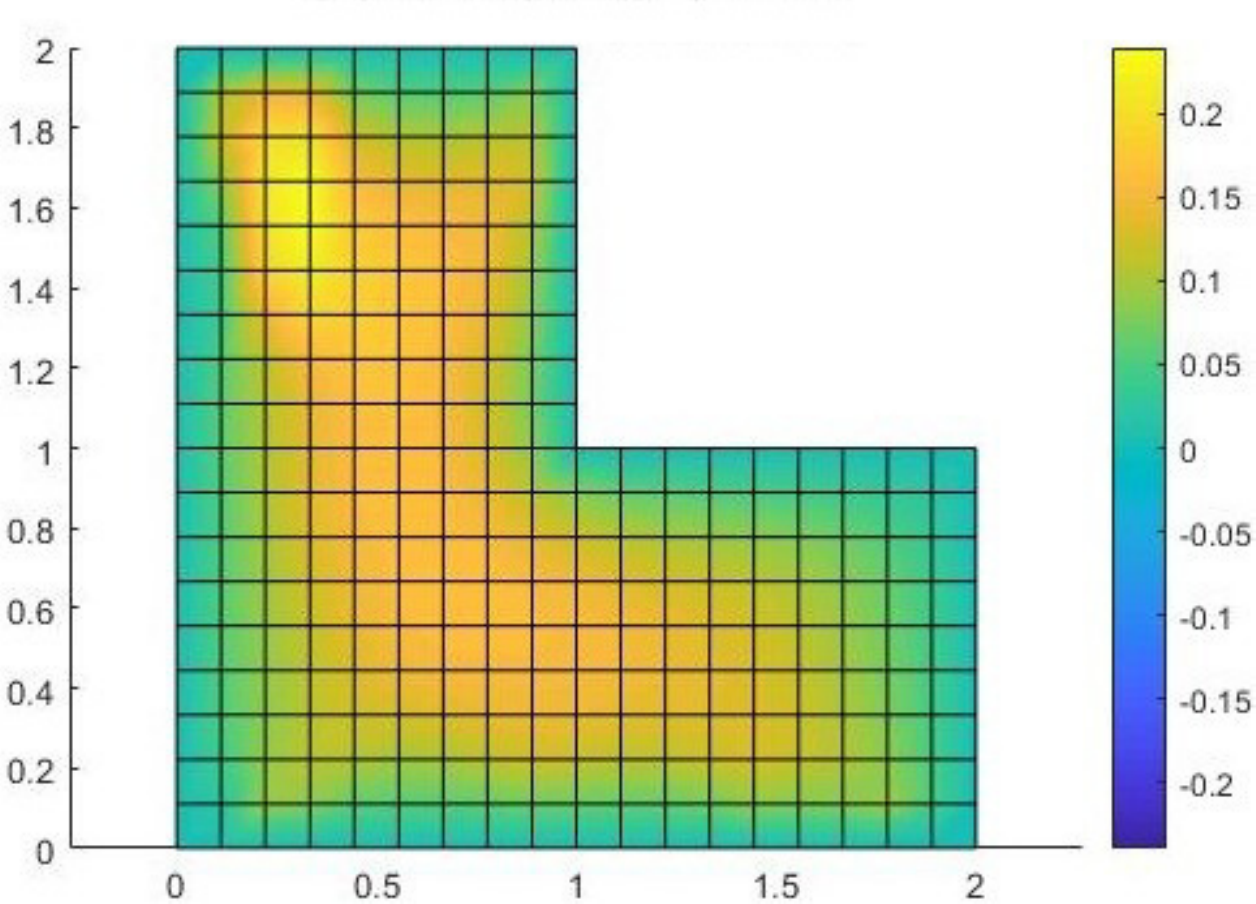}\qquad
  \includegraphics[width = 4.44cm]{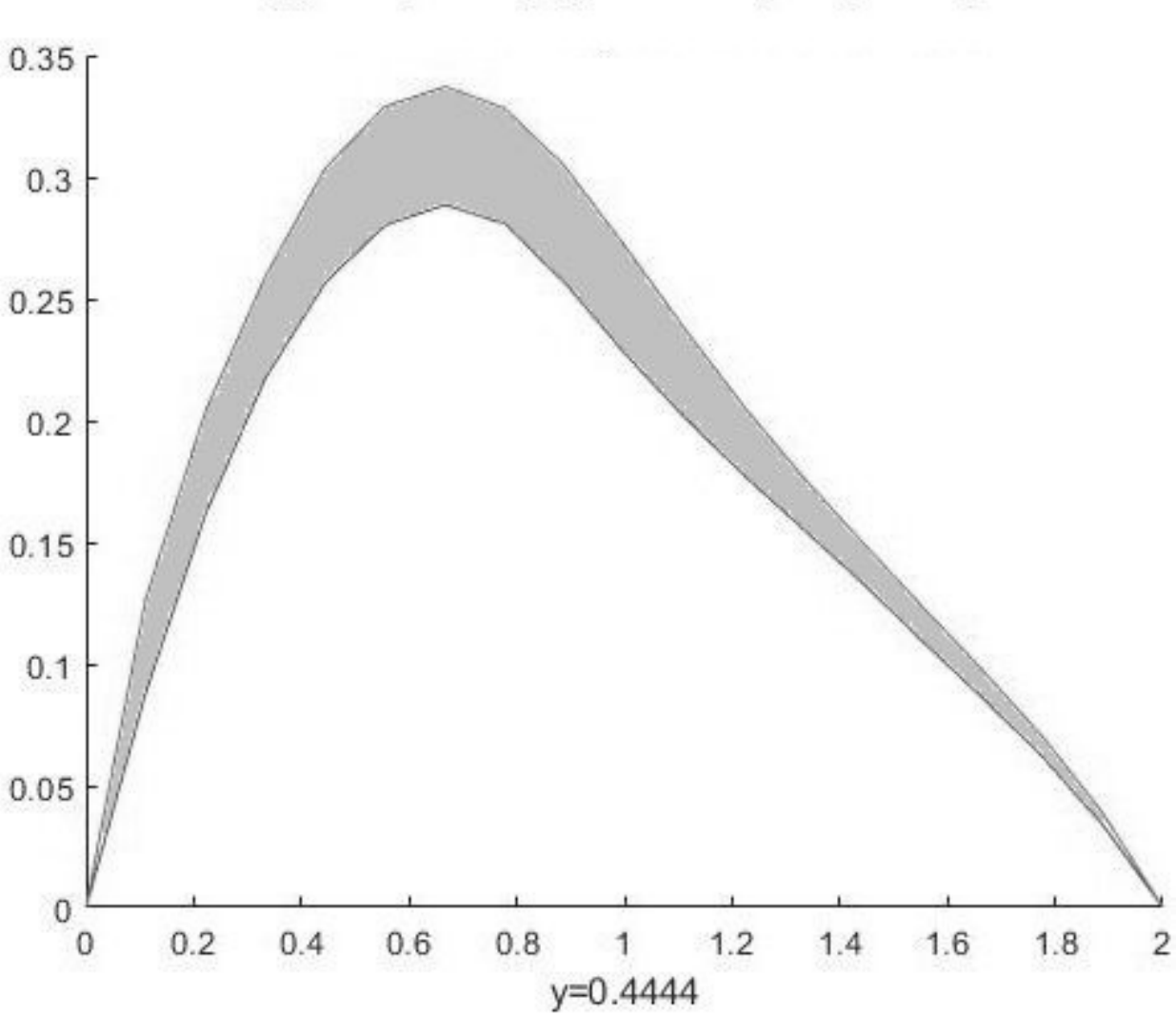}\qquad
  \raisebox{0pt}{\includegraphics[width = 4.44cm]{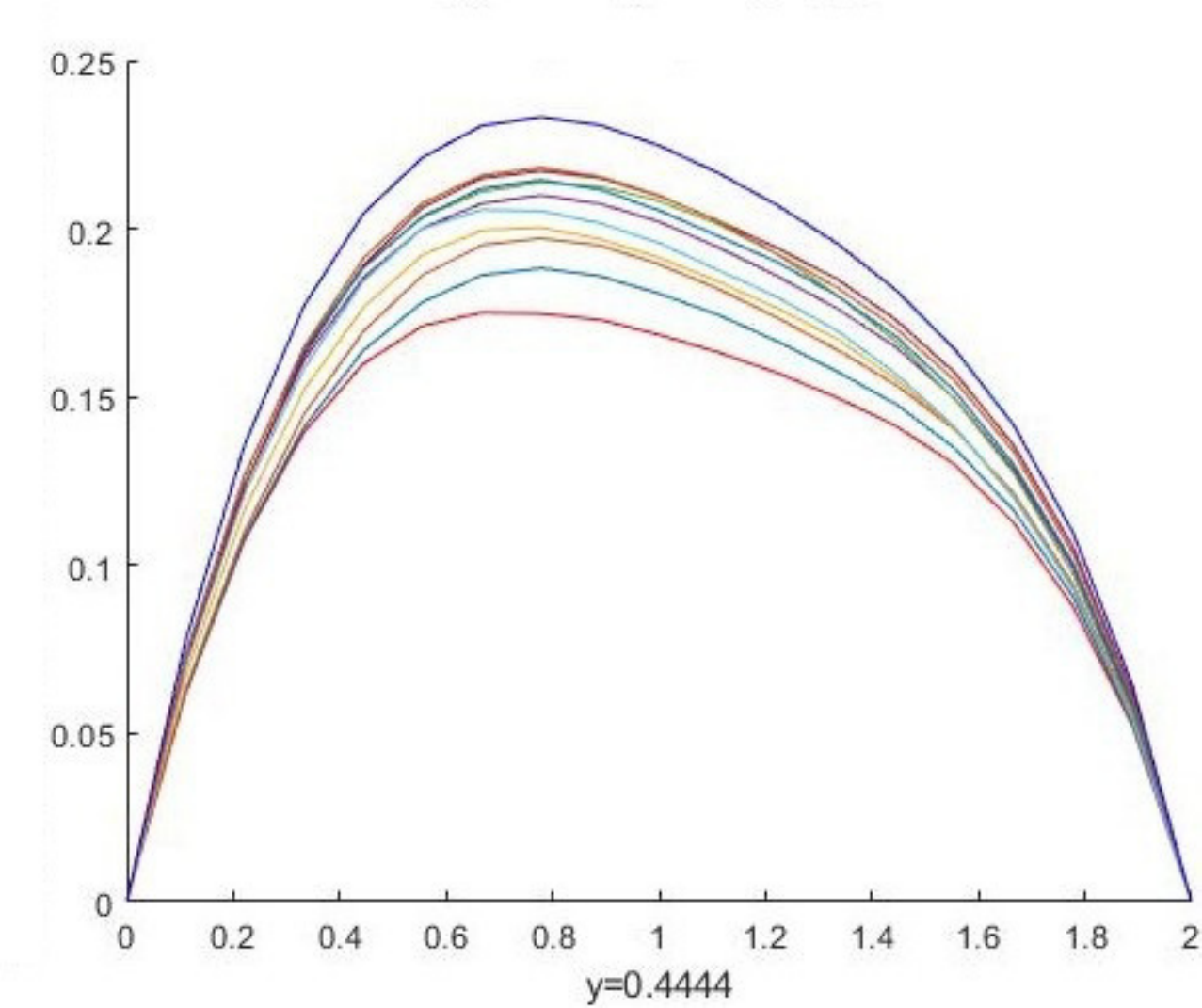}}
   \caption{Some aspects of the random set solution. Trajectory of displacement random field at fixed correlation length $\ell = 0.5$ (left);
   slice at $x_2 = 0.4444$ through a realization of the full random set solution (center); slice at $x_2 = 0.4444$ through the mean value interval field, showing individual mean fields for varying correlation lengths inside.}
  \label{fig:2}
\end{figure}

\emph{Numerical example.}
The purpose of this paragraph is to illustrate the method in a simple example, the displacement of an L-shaped membrane with non-uniform, random elastic properties.
This has the character of a text book example, and so all units are taken dimensionless. The L-shaped domain $D$ can be seen in Figure \ref{fig:2} (left). The equation for the displacement is (\ref{eq:ellprob}). A constant load $f(x)\equiv 1$ is assumed, and the random field $a(x)$ is taken of the form (\ref{eq:coeff}) with constant mean
$\mu(x)\equiv 1$, autocovariance function of the form (\ref{eq:ACF}), and field variance identically equal to one as well. In this elliptic problem, it suffices to guarantee that the realizations of the random field $a(x,\omega,\ell)$ stay above zero. Thus a lower cut-off function $\chi$ was chosen at the value 0.1. To produce the random set solution, the correlation length was taken from the interval $0.5\leq \ell\leq 1.5$.

Numerically, the L-shaped region was subdivided in a grid of size $18\times 18$, and problem (\ref{eq:ellprob}) was solved by the finite element method using a
{\sf Matlab}-program. To generate the random field, the Karhunen-Lo\`eve expansion truncated at $M = 130$ was used. The value of $a(x,\omega,\ell)$ entered in each finite element was chosen as the average of the values at its four nodes. A Monte Carlo sample size of 500 was used. At each realization of the underlying $2M$-dimensional Gaussian variable, see (\ref{eq:KLprocess}), the finite element program was executed with eleven values of $\ell$ between 0.5 and 1.5, stepsize 0.1.

As outlined above, this procedure results in an approximation to the random set solution $U(\omega)$. The solution is a random set in $H_0^1(D)$ and in $H^2(D')$ in any open subregion, actually everywhere except near the sharp inside corner \cite{1985:grisvard}. Consequently, each horizontal or vertical slice through the solution is a random set in the space of continuous functions. At fixed $x = (x_1,x_2)$, each realization of the solution is an interval. The mean (Aumann expectation) is an interval field. Figure \ref{fig:2} (left) shows a realization of the displacement random field at fixed correlation length $\ell = 0.5$. Realizations of the full random set solution can be visualized by plotting slices at fixed $x_1$ or $x_2$. Such a cut at $x_2 = 0.4444$ is shown in Figure \ref{fig:2} (center). Above each fixed $x_1$, the realization is an interval. Finally, the mean value field can also be depicted by slices. Figure \ref{fig:2} (right) shows a slice of this interval field at $x_2 = 0.4444$. The upper and lower bounding curves of the mean displacement are shown together with the mean solutions at correlation lengths 0.5, 0.6, \ldots, 1.5.

It is clear that the strategy outlined in the simple example (\ref{eq:ellprob}) can be extended to more general elliptic equations and systems, using the variational formulation (see \cite{1967:necas,1975:treves} for second order scalar elliptic equations and \cite{1976:duvaut} for more general elliptic systems of elasticity theory).

\section{Application to hyperbolic systems}
\label{sec:hyp}

Many equations from one-dimensional elastodynamics (and other fields, such as acoustics or hydrodynamics) can be written as a hyperbolic system
\begin{equation}\label{eq:hsys}
	\begin{array}{rcl}
	\left(\partial_t+a_i(x,t)\partial_x\right)u_i(x,t)&=&\sum_{j=1}^n f_{ij}(x,t)u_j(x,t) + g_i(x,t),\\
	u_i(x,0)&=&u_{0i}(x), \quad {\rm for}\ i = 1,\ldots, n,
	\end{array}	
\end{equation}
where $x\in\R$ and $t\in \R$ denote the spatial and time variable, respectively. The solution vector is $(u_1,\ldots,u_n)$, the driving term is $(g_1,\ldots,g_n)$, the coefficient matrix is $(f_{ij})_{i,j=1,\ldots,n}$.
The coefficients $a_i$ are assumed to be real-valued and signify the propagation speeds of the solution components.
The focus of this section will be on the case where the coefficients as well as the initial data are random fields/random sets on a probability space $(\Omega,\Sigma,P)$.

A typical example from elastodynamics concerns longitudinal waves in thin, non-homogeneous rods, described by the wave equation
\begin{equation}\label{eq:rods}
   \rho(x)\p_t^2u(x,t) - \p_x(E(x)\p_x u(x,t)) = q(x,t)
\end{equation}
where $u$ is the longitudinal displacement, $\rho$ the density, $E$ the modulus of elasticity, and $q$ the body force \cite[Chapter 2]{1975:graff}.
Equation (\ref{eq:rods}) is readily cast in the form of (\ref{eq:hsys})
\begin{equation}\label{eq:rodsys}
	\begin{array}{rcl}
	\left(\partial_t+a(x)\partial_x\right)u_1(x,t)&=& f(x,t)(u_2(x,t) - u_1(x,t)) + g(x,t),\\
	\left(\partial_t-a(x)\partial_x\right)u_2(x,t)&=& f(x,t)(u_2(x,t) - u_1(x,t)) + g(x,t)
	\end{array}	
\end{equation}
by setting $u_1(x,t) = (\p_t-a(x)\p_x)u(x,t)$, $u_2(x,t) = (\p_t+a(x)\p_x)u(x,t)$ with $a = \sqrt{E/\rho}$, $f = \p_x E/\rho - \p_xa/2$ and $g = q/\rho$.

A second example is provided by various scalar transport or advection-reaction equations in the form
\begin{equation}\label{eq:trans}
	\begin{array}{rcl}
	\p_t u(x,t) + a(x,t)\p_x u(x,t) &=& f(x,t)u(x,t) + g(x,t),\\
	u(x,0)&=&u_{0}(x).
	\end{array}	
\end{equation}
The construction of classical, deterministic solutions to the hyperbolic system (\ref{eq:hsys}) is a vector-valued version of the construction for the scalar equation (\ref{eq:trans}). Thus it will suffice to briefly indicate how
(\ref{eq:trans}) can be treated. The systems case can be found, e.g., in \cite[Section 13]{MOBook}. The solution can be constructed by the method of characteristics, which consists in solving the system of differential equations
\begin{equation}\label{eq:transchar}
\begin{array}{rcl}
	\frac{\partial}{\partial \tau}\gamma(x,t,\tau)&=&a\left(\gamma(x,t,\tau),\tau\right),\\
	\gamma(x,t,t)&=&x.
\end{array}
\end{equation}
This results in the characteristic curves $\tau \to \gamma(x,t,\tau)$, passing through the point $x$ at time $\tau = t$. Equation (\ref{eq:trans}) is equivalent to the integral equation
\begin{equation}\label{eq:transint}
	u(x,t)=u_0(\gamma(x,t,0))+\int_0^t\big(f(\gamma(x,t,\tau),\tau)u(\gamma(x,t,\tau),\tau) + g(\gamma(x,t,\tau),\tau)\big)\,d\tau,
\end{equation}
which in turn can be solved by successive iteration.

The issue to be addressed here is when $a$ and $f$, $g$ (and possibly $u_0$) are random fields or parametrized random fields on the probability space $(\Omega,\Sigma,P)$. In order to obtain the solution as a random field or as a random set, additional properties of the coefficients are required.

\emph{The random field case.}
The starting observation is that equation (\ref{eq:transchar}) turns into the random differential equation
\[
  \frac{\partial}{\partial \tau}\gamma(x,t,\tau,\omega) = a\left(\gamma(x,t,\tau,\omega),\tau,\omega\right),\quad \gamma(x,t,t,\omega)=x.
\]
Following the theory of random differential equations \cite{1972:bunke,2017:han}, a continuous random field solution $\gamma(x,t,\tau,\omega)$ can be obtained provided the transport coefficient $a(x,t,\omega)$ is pathwise Lipschitz continuous and globally bounded.
The next issue is the composition of the random fields $f$, $g$, $u_0$ with $\gamma$, for example, $(x,t,\tau,\omega)\to f(\gamma(x,t,\tau,\omega),\tau,\omega)$. This composition is jointly measurable (in all four variables) and pathwise continuous, if both $\gamma$ and $f$ are pathwise continuous random fields. Finally, the integral equation (\ref{eq:transint}), now with random field input, has to be solved by successive approximation (on suitable bounded domains of determinacy), taking account of the measurability at each iteration step. The precise requirements are:
\begin{itemize}
\item For fixed $\omega\in\Omega$, $a(x,t,\omega)$ is locally Lipschitz in $(x,t)$;
\item $|a(x,t,\omega)|$ has a uniform global bound independently of $(x,t,\omega)$;
\item $|f(x,t,\omega)|$ is locally bounded in $(x,t)$ independently of $\omega$;
\item $f$, $g$ and $u_0$ are pathwise continuous random fields.
\end{itemize}
Under these assumptions, the integral equation (\ref{eq:transint}) has a pathwise continuous random field solution $u(x,t,\omega)$. Note that the assumptions on $a$ and $f$ can be enforced by a cut-off procedure similar to (\ref{eq:coeff}).
If, in addition, the random fields $a$, $f$, $g$, $u_0$ are pathwise differentiable, $u(x,t,\omega)$ is a pathwise classical solution to the differential equation (\ref{eq:trans}). When solving, e.g., the wave equation (\ref{eq:rods}), one has to take care that the fields $\rho$ and $E$ are sufficiently regular so that the coefficients in system (\ref{eq:rodsys}) have the required properties.

Another approach to hyperbolic systems with random field coefficients, based on weak solutions, can be found in \cite{2016:mishra}. It applies in any space dimension and leads to a different set of requirements on the coefficients and data.

\emph{The random set case.}
Assume that some or all of the random fields $a$, $f$, $g$, $u_0$ depend on interval parameters. As in the previous sections, this leads to a set-valued solution $\omega\to U(\omega)$ to (\ref{eq:trans}) in a suitable space to be determined.
In order to prove that the set-valued solution is a random set, one may follow exactly the same strategy as in Section~\ref{sec:elastostatics}.
Only the continuous dependence of the solution on the coefficients and data requires different arguments, as outlined next. The continuity of the map $(a,f,g,u_0) \to u = u_{(a,f,g,u_0)}$ can be split up into the continuity of the maps
\[
   a\to\gamma\quad {\rm and}\quad (\gamma,f,g,u_0)\to u.
\]
To be specific, start with a bounded and closed initial interval $K_0 = [-\kappa,\kappa]$ and a time horizon $T>0$. Denote by $c$ the uniform global bound of $|a(x,t,\omega)|$, as listed under the requirements. The region
\[
  K_T = \{(x,t)\in\R^2: -T\leq t \leq T, |x|\leq\kappa - c|t|\}
\]
is a domain of determinacy for $K_0$. That is, when $(x,t)\in K_T$, all characteristic curves $\tau\to\gamma(x,t,\tau,\omega)$ joining $(x,t)$ with the $x$-axis lie also in $K_T$. Denote by ${\rm Lip}(K_T)$ the space of Lipschitz continuous functions on $K_T$ and by $\Gamma(K_T\times[-T,T])$ the subset of $\cC(K_T\times[-T,T])$ which comprises all characteristic curves $\gamma$ arising from elements $a\in{\rm Lip}(K_T)$. The first step is to show that the map
\[
   {\rm Lip}(K_T)\to \Gamma(K_T\times[-T,T]): a\to\gamma
\]
is continuous, when the spaces are equipped with the supremum norm. The second step is to prove that the map
\[
  \Gamma(K_T\times[-T,T])\to \cC(K_T): (\gamma,f,g,u_0)\to u
\]
is continuous in the supremum norm. The latter is a rather standard argument in the theory of differential equations using the Gronwall lemma. Details are elaborated in \cite{2020:nedeljkovic}.  Since the initial choice of $K_0$ was arbitrary, the continuity of the map $(a,f,g,u_0) \to u = u_{(a,f,g,u_0)}$ with values in $\cC(\R^2)$ follows. Thus if the random fields $(a,f,g,u_0)$ depend continuously on a (possibly multidimensional) parameter $\lambda\in\Lambda$, the same arguments as in Sections~\ref{sec:parrf} and \ref{sec:elastostatics} show that
\[
   U(\omega) = \{u_{(a(\omega,\lambda),f(\omega,\lambda),g(\omega,\lambda),u_0(\omega,\lambda))}:\lambda\in\Lambda\}
\]
is a random set in $\cC(\R^2)$.

\emph{The multidimensional case:}
The generalization to the transport equation in higher space dimensions,
\[
   \partial_t u(x,t) +a(x,t)\cdot {\rm grad\,}u(x,t) = f(x,t)u(x,t) + g(x,t)
\]
is straightforward, using the method of characteristics. Here $x\in\R^d$, $a = (a_1,\ldots, a_d)$, $u$, $f$ and $g$ are still scalar functions. In the multidimensional analogues of (\ref{eq:hsys}), (\ref{eq:rods}), it is more challenging to establish continuous dependence of the solution on the data. One approach could be to formulate the problem in terms of semigroup theory and then invoke Kato's perturbation results \cite{1976:kato}. For the system of linear acoustics this has been done, e.g., in \cite[Section 6]{2015:garetto}. Another approach  takes recourse to the Fourier integral operator representation of the solution \cite{2020:MOSchwarz}.

\section{Further results, numerical aspects}
\label{sec:further}

\subsection{The parametric point of view}
\label{subsec:parametric}

Given a parametrized family $a_\lambda, \lambda\in\Lambda$, of random variables or random fields, and a map $a_\lambda\to u_{a_\lambda}$ describing the response of a system, the approach favored in this paper has been to first construct the random set
\begin{equation}\label{eq:randomsetcollected}
   U(\omega) = \{u_{a_\lambda(\omega)}:\lambda\in\Lambda\}
\end{equation}
and then to define lower and upper probabilities of events $B$ by
\[
  \underline{P}(B) = P(U_-(B)),\quad \overline{P}(B) = P(U^-(B)).
\]
An alternative approach, the \emph{parametric viewpoint}, would consist in skipping the random set step and defining lower and upper probabilities induced by the family $a_\lambda, \lambda\in\Lambda$, directly:
\[
  P_{\rm low}(B) = \inf_{\lambda\in\Lambda} P(u_{a_\lambda}(B)),\quad P_{\rm upp}(B) = \sup_{\lambda\in\Lambda}P(u_{a_\lambda}(B)).
\]
This approach gives tighter bounds:
\[
  \underline{P}(B) \leq  P_{\rm low}(B) \leq P_{\rm upp}(B) \leq  \overline{P}(B),
\]
but deprives one of using the concepts of the theory of random sets. The numerical computation of the two types of upper and lower probabilities are different, but of comparable effort. For a comparison of the two approaches, from the conceptual and numerical points of view, see \cite{2017:alvarez,2016:fetz}.

\subsection{Numerical aspects}
\label{subsec:num}

As is well-known, the computational effort when implementing heterogeneous uncertainty models is high: it involves both Monte Carlo sampling and the computation of interval bounds. Suppose the quantity of interest, given through  a model response function $a_\lambda \to u_{a_\lambda}$ as above, is one-dimensional. This could be the model response at a single point $x$ or the value of a performance criterion. Further, suppose the lower and upper distribution functions $\underline{F}(b)$, $\overline{F}(b)$ of this quantity are sought.

The procedure in the random set framework is as follows: If the hypotheses of the previous sections are met, the quantity of interest (\ref{eq:randomsetcollected}) is a random interval $U(\omega) =  [\underline{u}(\omega),\overline{u}(\omega)]$, and the lower and upper distribution functions are simply given by
$\underline{F}(b) = P(\overline{u}\leq b)$, $\overline{F}(b) = P(\underline{u}\leq b)$.
The following algorithm produces a numerical approximation to these lower and upper distribution functions. From there, further probabilistic indicators such as bounds on exceedence probabilities or moment bounds can be computed.
It is assumed that the parameter set $\Lambda$ indexing the family of random variables $a_\lambda$ is a $d$-dimensional interval (or a sufficiently regular subset of such an interval). It is also assumed that the random variables $a_\lambda$ are defined on a common probability space $(\Omega,\Sigma,P)$.
The following steps are required:
\begin{itemize}
\item Discretize $\Lambda$ in a grid $\lambda_1,\ldots,\lambda_M$.

\item Generate a sample $\omega_1,\ldots, \omega_N$ distributed according to the underlying probability $P$ on the common probability space $\Omega$.

\item For each $\omega_k$, estimate
\[
\underline{u}(\omega_k) = \min_{i=1,\ldots,M}u_{a_{\lambda_i}}(\omega_k),\quad
   \overline{u}(\omega_k) = \max_{i=1,\ldots,M}u_{a_{\lambda_i}}(\omega_k).
\]
Each $[\underline{u}(\omega_k),\overline{u}(\omega_k)]$ is a realization of the random interval $[\underline{u},\overline{u}]$.

\item Then $\underline{u}(\omega_k)$ and $\overline{u}(\omega_k)$, $k = 1,\ldots, N,$ are samples of $\underline{u}$ and $\overline{u}$, respectively.
Their empirical cumulative distribution functions provide estimates for $\underline{F}(b)$ and $\overline{F}(b)$.
\end{itemize}
In the parametric point of view, the assumptions on the index set $\Lambda$ are the same. It is not necessary to require that the random variables $a_\lambda$ are defined on the same probability space. Denote by $P_\lambda$ the probability distribution of $a_\lambda$.
The algorithm is as follows:
\begin{itemize}
\item Discretize $\Lambda$ in a grid $\lambda_1,\ldots,\lambda_M$.

\item For each $\lambda_i$, generate a sample ${a_{\lambda_i}}(1),\ldots, {a_{\lambda_i}}(N)$ of the random variable ${a_{\lambda_i}}$ (that is, generate a sample of size $N$ from the distribution $P_{\lambda_i}$).

\item Estimate the probability $P_{\lambda_i}((-\infty,b])$ using the sample, for example by means of the Monte Carlo estimator $\frac1{N}\sum_{j=1}^N \chi_{(-\infty,b]}(a_{\lambda_i}(j))$ where $\chi_{(-\infty,b]}$ denotes the characteristic function of the interval $(-\infty,b]$.

\item Estimate
\[
F_{\rm low}(b) \approx \min_{i=1,\ldots,M}P_{\lambda_i}((-\infty,b]),\quad F_{\rm upp}(b) \approx \max_{i=1,\ldots,M}P_{\lambda_i}((-\infty,b]).
\]
\end{itemize}
In these crude forms, both algorithms require $M\cdot N$ evaluations of the (generally expensive) model function $u$. Accordingly, many approaches have been devised to reduce computational cost. In the random set case, a stochastic response surface, like a polynomial chaos expansion \cite{1991:ghanem}, reduces the number $M$ to a level required for the accuracy of the response surface. Evaluating a Monte Carlo sample of size $N$, given the response surface, is cheap \cite{2015:oberguggenberger}. For further sampling concepts for random sets, see \cite{2006:alvarez}.

Here is a short list of methods currently being developed for the reduction of computational cost in imprecise probability models, in particular, for imprecise random fields:
Approximations by \emph{inner and outer bounds} \cite{2009:alvarez,2017:alvarez,2018:alvarez,2004:tonon,2020:zhang}; propagation methods for \emph{$p$-boxes} \cite{2021:FaesBroggi,2021:FaesDaub,2022:gray,2020:sadeghi,2019a:wei,2019b:wei}; methods employing \emph{polynomial chaos expansion} \cite{2020:liu,2014:oberguggenberger,2017:schoebi}; \emph{probability plots}, an enhancement of the first order reliability method (FORM) proposed by \cite{2017:hurtado}.
A very efficient method appears to be \emph{advanced line sampling} \cite{2015:deangelis,2020:song} that intertwines the two required loops and can reduce $M$ and $N$ simultaneously. In the parametric point of view, a promising method has become \emph{importance sampling}, or rather \emph{reweighting} a single importance Monte Carlo sample, which means $M = 1$. Still $N$ model evaluations are required \cite{2017:fetz,2019:fetz,2017:troffaes,2018:troffaes,2018:troffaesfetzMO,2018:zhang}. This can of course be aided by a deterministic response surface.

\subsection{Modelling correlations}
\label{subsec:corr}

The question of how to model correlations in hybrid probability-interval uncertainty analysis has occurred implicitly at various places in this article. The purpose of this subsection is to spell out more explicitly what kind of approaches are possible. Essentially, there are three ways to formulate correlations.

\emph{Correlations introduced by the underlying probabilistic model.} In the approach to random sets based on families of random variables, correlations can be defined through their joint distributions. These correlations are carried along when replacing distributional hyperparameters by intervals. This is generally the case in random field models with interval correlation lengths as e.g in \cite{2015:oberguggenberger} as well as in Section~\ref{sec:elastostatics}. More generally, random fields with families of correlation structures have been considered in \cite{2019:faes}, especially analyzing the corresponding Karhunen-Lo\`eve expansions with interval valued eigenvalues and eigenfunctions. In models with interval parameters and random field excitation, correlations are introduced through the random field, but become interval valued in output and/or state variables. Such a situation has been investigated in connection with spectral finite elements (finite elements combined with polynomial chaos expansions) in \cite{2016:do}.

\emph{Copulas and multivariate capacities.} A very interesting and practical approach is due to \cite{2006:alvarez,2009:alvarez}, who introduced the notion of \emph{random sets of indexable type}. These random sets are defined on $\Omega = (0,1]^d$ and have values in a target product space $\cA = \cA_1\times\cdots\times \cA_d$. The focal elements are of product form $A_1(\omega_1)\times \cdots \times A_d(\omega_d)$ where $\omega = (\omega_1,\ldots,\omega_d)\in \Omega$ and each $A_j(\omega_j)$ belongs to a suitable family of subsets of $\cA_j$ (closed subsets for the purpose of the present paper). A correlation structure can be imposed on a random set of indexable type by prescribing the joint probability measure on $\Omega$ with the help of a $d$-dimensional copula. The approach has been further developed in connection with simulation methods in \cite{2018:alvarez}.

Alternatively, copulas can also be introduced so as to act on the target space, rather than on the underlying probability space $\Omega$. A bit of terminology is in order. The upper probability of a Borel subset $B\subset \cA$ defined by (\ref{eq:uppprob}) can be viewed as a so-called \emph{capacity functional} \cite{2005:molchanov}, say on the family $\cK(\cA)$ of compact subsets of $\cA$, while the lower probability defines a \emph{containment functional}. In the multidimensional case, there are two possible ways of extending capacity functionals, following \cite{2019:schmelzer}. The first one is to define a capacity functional on the family of sets belonging to $\cK(\cA_1)\times \cdots \times \cK(\cA_d)$ directly, in which case it is termed a \emph{multivariate capacity functional}. The second choice is to define a capacity functional on the larger family of sets $\cK(\cA_1\times\cdots\times\cA_d)$, referred to as a \emph{joint capacity functional}. It has been shown in \cite{2015:schmelzer} that multivariate containment functionals can be reconstructed from their univariate marginals by means of a family of copulas (but not by a single copula, in general). On the other hand, joint capacity functionals can always be reconstructed from their marginal capacities by means of a single copula, provided the marginal capacities are maxitive \cite{2019:schmelzer} (this is the case if the marginals arise from normalized fuzzy sets as indicated in Subsection~\ref{subsec:rs}). The relation of this approach to the one of \cite{2009:alvarez} and others has been discussed in \cite{2015a:schmelzer}.

\emph{Structured focal elements.} The third possibility to model dependence is through the structure of the focal elements. The geometric approach would replace, say in the case of two variables, the rectangles $A_1(\omega_1)\times A_2(\omega_2)$ defining the joint focal sets by rhombic shapes \cite{2010:bernardini:book,2010:tonon}. Alternatively, any random set delivers a family of probability distributions, see e.g. \cite{2013:oberguggenberger}. One may keep the focal elements as rectangles but restrict the admissible joint probability distributions \cite{2004:fetz}. In the case of interval fields, one may use explicitly defined basis functions multiplied by interval factors \cite{2013:verhaeghe} or \emph{spatial dependency functions} \cite{2015:sofi}, which are the counterpart to Karhunen-Lo\`eve expansions in the interval setting.

\section{Summary and conclusions}
\label{sec:conclusion}
The paper addressed mathematical and modelling issues when heterogeneous uncertainty is entered in engineering models. The model inputs were assumed to consist of random variables or random fields, combined with interval parameters. This led to model outputs which were both random and set-valued. This set-up is not new; it actually has been widely adopted in the engineering literature \cite{2016:do,2010:gao,2018:lue,2016:muscolino,2015:oberguggenberger,2008:MOFellin,2010:schmelzer:adam,2022:Sofi,2020:Sofi,2020:song}, see also the literature review in the introduction. However, few papers have addressed the question whether the set-valued model output satisfies the measurability conditions which make it a random set in its proper sense \cite{2017:alvarez,2019:wurzer,2010:schmelzer,2013:schmelzer}. The main contribution of this paper was to show how the theoretical measurability questions can be addressed, and that the continuous dependence of the solution and the involved random and interval quantities on their parameters is the essential ingredient. This is important in as much as computing upper and lower probabilities in the sense of (\ref{eq:uppprob}), (\ref{eq:lowprob}) is only justified when the model output is a random set.

It was assumed here that the underlying probability spaces were infinite and that the quantity of interest was the model output itself. If one is only interested in upper and lower bounds on certain probabilities or moments, one may avoid collecting the output in a random set. Rather, one may work directly with the families of random variables defining the input, propagating them through the model and computing the desired probability bounds. In this parametric approach, which has been adopted and discussed  e.g. in \cite{2017:alvarez,2019:faes,2017:fetz,2018:troffaes}, the measurability issue is much simpler: it reduces to proving that the propagated variables are measurable (that is, still random variables), which is guaranteed if the input-output map is continuous.

Another simple situation in which the measurability question does not arise is the case of finite random sets, because on a finite probability space (with the power set as the $\sigma$-algebra) every map is measurable.

The focus was on model outputs that were generated through solutions to elliptic and hyperbolic partial differential equations. The most difficult part was the proof that the model output depended continuously on the input parameters and hyperparameters. There is no general method of proof; different arguments are needed for each type of equations. It is expected that new arguments will be needed when considering other types of equations, such as parabolic equations. As the paper was oriented toward the theoretical background of a widely adopted framework, only a simple illustrative example was given. Many larger scale applications can be found in the papers quoted above. For the same reason, computational aspects have only been briefly addressed. For more extensive analyses of the computational effort, the reader is referred to \cite{2018:alvarez,2021:FaesDaub,2017:fetz,2016:fetz}.

An important future research question is the modelling of correlations. It is expected that the framework of \cite{2019:schmelzer} can encompass most other approaches. In particular, it will be interesting to investigate how it applies in the case of interval fields.

\section*{Acknowledgements}
This work was supported by the FWF doctoral program Computational Interdisciplinary Modelling W1227 of the Austrian Science Fund as well as a doctoral scholarship of the University of Innsbruck. Thanks are due to Thomas Fetz for providing the {\sf Matlab}-code solving the deterministic FE-problem (\ref{eq:ellprob}).
The authors are grateful to the two referees for valuable comments. Thanks to the first referee, the important issue of correlation modelling has been included. The detailed and critical comments of the second referee led to many clarifications and to substantial improvements of the presentation.

\section*{Appendix: Explicit construction of infinite probability spaces}

In Subsection~\ref{subsec:contrf} it was necessary to introduce common probability spaces for random fields with different variances and correlation lengths, first for random fields defined through the Langevin equation, and second for the infinitely many Gaussian random variables arising in the Karhunen-Lo\`eve expansion. The appendix briefly explains the mathematical details.

\emph{Wiener space $(\Omega_W,\Sigma_W,P_W)$.} Here $\Omega_W = \cC([0,\infty))$ is the space of continuous functions $\omega = \omega(x)$ on the half-line. It is a metric space; the metric is generated by the semi-norms
$
  \|\omega\|_{a,b} = \max_{a\leq x \leq b}|\omega(x)|,
$
$0\leq a < b < \infty$. Then $\Sigma_W$ is the corresponding Borel $\sigma$-algebra. The probability measure $P_W$ is the measure on $\Omega_W$ which is generated by the transition probabilities of Wiener process \cite{1987:stroock}.
In this setting, Wiener process $W:[0,\infty)\times \Omega_W \to\R$ is defined through the point evaluations $W(x,\omega) = \omega(x)$. Note that the realizations of the process are identified with the random elements, as is usual practice when constructing probability spaces on which stochastic processes live. The one-point and two-point probabilities are, for points $0 < x < y$,
\[
   P_W(a\leq \omega(x)\leq b) = \int_{a}^{b} p(x,z)\dd z\quad \mbox{with}\quad p(x,z) = \frac{1}{\sqrt{2\pi x}}\exp\big(-\frac{z^2}{2x}\big)
\]
\[
   P_W(a\leq \omega(x) \leq b, c \leq \omega(y) \leq d) = \int_a^b p(x,z_1)\int_c^d p(y-x,z_2 - z_1)\dd z_2\dd z_1
\]
and similarly for the joint probabilities at finitely many points.

\emph{Infinite products of Gaussian spaces.} For $k = 1,2,3,\ldots$, let $\Omega_k = \R$ with the Borel $\sigma$-algebra $\Sigma_k$ and the standard Gaussian measure $P_k$,
\[
P_k(B_k) = \frac1{\sqrt{2\pi}}\int_{B_k} \exp(-\frac12 \omega_k^2)\dd \omega_k
\]
for any Borel subset $B_k\subset\R$. The spaces $(\Omega_k,\Sigma_k,P_k)$ are all identical copies of standard Gaussian space, but the index is needed to form their product.
The standard Gaussian variable on $\Omega_k = \R$ is defined through the identity map ${\rm id}(\omega_k) = \omega_k$. Again, the realizations of the random variable are identified with the random elements.

Let $\Omega$ be the infinite product of the spaces $\Omega_k$. The random elements are sequences $\omega = (\omega_1,\omega_2,\omega_3,\ldots)$. The probability measure on $\Omega$ is the product measure
\[
   P(B_1\times B_2 \times \ldots \times B_n) = \prod_{k=1}^n P_k(B_k) = (2\pi)^{-n/2} \prod_{k=1}^n \int_{B_k} \exp(-\frac12 \omega_k^2)\dd \omega_k
\]
for any finite collection of Borel sets $B_1, B_2,\ldots, B_n$.
The infinite dimensional random variable $\xi:\Omega\to\R$, $\xi(\omega) = (\xi_1(\omega),\xi_2(\omega),\ldots)$ is given through the projections to $\Omega_k$, that is,
$\xi_k(\omega) = \xi_k(\omega_1,\omega_2,\ldots) = \omega_k$.


\begin{thebibliography}{}

\bibitem{2006:alvarez}
D.A. Alvarez,
\newblock {On the calculation of the bounds of probability of events using infinite random sets},
\newblock {International Journal of Approximate Reasoning}, 43 (2006) 241--267.
https://doi.org/10.1016/j.ijar.2006.04.005.

\bibitem{2009:alvarez}
D.A. Alvarez,
\newblock {A Monte Carlo-based method for the estimation of lower and upper probabilities of events
using infinite random sets of indexable type},
\newblock {Fuzzy Sets and Systems}, 160 (2009) 384--401.
https://doi.org/10.1016/j.fss.2008.08.006.

\bibitem{2017:alvarez}
D.A. Alvarez, J.E. Hurtado, J. Ram\'{i}rez,
\newblock {Tighter bounds on the probability of failure than those provided by random set theory},
\newblock {Computers \& Structures}, 189 (2017) 101--113.
http://dx.doi.org/10.1016/j.compstruc.2017.04.006.

\bibitem{2018:alvarez}
D.A. Alvarez, F. Uribe, J.E. Hurtado,
\newblock {Estimation of the lower and upper bounds on the probability of failure using subset simulation and random set theory},
\newblock {Mechanical Systems and Signal Processing}, 100 (2018) 782--801.
https://doi.org/10.1016/j.ymssp.2017.07.040.

\bibitem{1974:arnold}
L. Arnold,
\newblock {Stochastic Differential Equations: Theory and Applications},
\newblock Wiley-Interscience, New York-London-Sydney, 1974.

\bibitem{2013:beer:ferson:kreinovich}
M. Beer, S. Ferson, V. Kreinovich,
\newblock Imprecise probabilities in engineering analysis,
\newblock {Mechanical Systems and Signal Processing}, 37 (2013) 4--29.
http://dx.doi.org/10.1016/j.ymssp.2013.01.024.

\bibitem{2010:bernardini:book}
A. Bernardini, F.~Tonon,
\newblock {Bounding Uncertainty in Civil Engineering -- Theoretical  Background},
\newblock Springer, Berlin, 2010.

\bibitem{1972:bunke}
H. Bunke,
\newblock {Gew\"{o}hnliche {D}ifferentialgleichungen mit zuf\"{a}lligen {P}arametern},
\newblock Akademie-Verlag, Berlin, 1972.

\bibitem{1977:castaing}
C. Castaing, M. Valadier,
\newblock {Convex analysis and measurable multifunctions},
\newblock Springer-Verlag, Berlin-New York, 1977.

\bibitem{2015:deangelis}
M. de Angelis, E. Patelli, M. Beer,
\newblock {Advanced Line Sampling for efficient robust reliability analysis},
\newblock {Structural Safety}, 52 (2015) 170--182.
http://dx.doi.org/10.1016/j.strusafe.2014.10.002.

\bibitem{2016:do}
M.D. Do, W. Gao, C. Song, M. Beer,
\newblock {Interval spectral stochastic finite element analysis of structures
with aggregation of random field and bounded parameters},
\newblock {International Journal for Numerical Methods in Engineering}, 108 (2016) 1198--1229.
https://doi.org/10.1002/nme.5251.

\bibitem{1976:duvaut}
G. Duvaut, J.-L. Lions (1976),
\newblock {Inequalities in Mechanics and Physics},
\newblock Springer-Verlag, Berlin-New York, 1976.

\bibitem{2021:FaesBroggi}
M. Faes, M. Broggi, K.-K. Phoon, M. Beer,
\newblock Distribution-free {P}-box processes: definition and simulation,
\newblock in: A. Sofi, G. Muscolino, R.L. Muhanna (Eds.), {REC 2021: Proceedings of the 9th International Workshop on Reliable Engineering Computing}, 2021,
pp. 380--394. http://ww2new.unime.it/REC2021/proceedings/REC2021\_Proceedings.pdf.

\bibitem{2021:FaesDaub}
M. Faes, M. Daub, S. Marelli, E. Patelli, M. Beer,
\newblock {Engineering analysis with probability boxes: A review on computational methods},
\newblock {Structural Safety}, 93 (2021) 102092.
https://doi.org/10.1016/j.strusafe.2021.102092.

\bibitem{2017:faes}
M. Faes, D. Moens,
\newblock {Identification and quantification of spatial interval uncertainty in numerical models},
\newblock {Computers and Structures}, 192 (2017) 16--33.
http://dx.doi.org/10.1016/j.compstruc.2017.07.006.

\bibitem{2019:faes}
M. Faes, D. Moens,
\newblock {Imprecise random field analysis with parametrized kernel functions},
\newblock {Mechanical Systems and Signal Processing}, 134 (2019) 106334.
https://doi.org/10.1016/j.ymssp.2019.106334.

\bibitem{2003:ferson}
S. Ferson, L. Ginzburg, V. Kreinovich, D. Myers, K. Sentz,
\newblock Constructing probability boxes and {D}empster-{S}hafer structures,
\newblock Technical Report SAND2002-4015,
Sandia National Laboratories, Albuqueque and Livermore, 2003.
https://doi.org/10.2172/809606.

\bibitem{2017:fetz}
T. Fetz,
\newblock Efficient Computation of upper probabilities of failure,
\newblock in: C. Bucher, B.R. Ellingwood, D.M. Frangopol (Eds.), {Proceedings of the 12th International Conference on Structural Safety and Reliability (ICOSSAR 2017)},
TU-MV Media Verlag GmbH, Wien, 2017, pp. 493--502.
https://owncloud.tuwien.ac.at/index.php/s/DnQ2sWRcnhT8Zfq\#pdfviewer.

\bibitem{2019:fetz}
T. Fetz,
\newblock {Improving the convergence of iterative importance sampling for computing upper and lower expectations},
\newblock {Proceedings of Machine Learning Research}, 103 (2019) 185--193.
https://proceedings.mlr.press/v103/fetz19a.html.

\bibitem{2004:fetz}
T. Fetz, M. Oberguggenberger,
\newblock {Propagation of uncertainty through multivariate functions in the framework of sets of probability measures},
\newblock {Reliability Engineering \& System Safety}, 85 (2004) 73--87.
https://doi.org/10.1016/j.ress.2004.03.004.

\bibitem{2016:fetz}
T. Fetz, M. Oberguggenberger,
\newblock {Imprecise random variables, random sets, and Monte Carlo simulation},
\newblock {International Journal of Approximate Reasoning}, 78 (2016) 252--264.
https://doi.org/10.1016/j.ijar.2016.06.012.

\bibitem{2010:gao}
W. Gao, C. Song, F. Tin-Loi,
\newblock {Probabilistic interval analysis for structures with uncertainty},
\newblock {Structural Safety}, 32 (2010) 191--199.
https://doi.org/10.1016/j.strusafe.2010.01.002.

\bibitem{2015:garetto}
C. Garetto, M. Oberguggenberger,
\newblock {Symmetrisers and generalised solutions for strictly hyperbolic systems with singular coefficients},
\newblock {Mathematische Nachrichten}, 288 (2015) 185--205.
https://doi.org/10.1002/mana.201400192.

\bibitem{1991:ghanem}
R. Ghanem, P. Spanos,
\newblock {Stochastic {F}inite {E}lements: {A} {S}pectral {A}pproach},
\newblock Springer, New York, 1991.

\bibitem{2002:goodman}
I.R. Goodman, H.T. Nguyen,
\newblock Fuzziness and randomness,
\newblock in: C. Bertoluzza, M. Gil, D.A. Ralescu (Eds.), {Statistical Modeling, Analysis and Management of Fuzzy Data},
Physica-Verlag, Heidelberg, 2002, pp. 3--21.

\bibitem{1975:graff}
K.F. Graff,
\newblock {Wave Motion in Elastic Solids},
\newblock Oxford University Press, Oxford, 1975.

\bibitem{2022:gray}
A. Gray, A. Wimbush, M. {de Angelis}, P.O. Hristov, D. Calleja, E. Miralles-Dolz, R. Rocchetta,
\newblock {From inference to design: A comprehensive framework for uncertainty quantification in engineering with limited information},
\newblock {Mechanical Systems and Signal Processing}, 165 (2022) 108210.
https://doi.org/10.1016/j.ymssp.2021.108210.

\bibitem{2010:grigoriu}
M. Grigoriu,
\newblock {Probabilistic models for stochastic elliptic partial differential equations},
\newblock {Journal of Computational Physics}, 229 (2010) 8406--8429.
https://doi:10.1016/j.jcp.2010.07.023.

\bibitem{1985:grisvard}
P. Grisvard,
\newblock {Elliptic Problems in Nonsmooth Domains}.
\newblock Pitman Advanced Publishing Program, Boston, 1985.

\bibitem{2017:han}
X. Han, P.E. Kloeden,
\newblock {Random Ordinary Differential Equations and Their Numerical Solution},
\newblock Springer, Singapore, 2017.

\bibitem{2017:hurtado}
J.E. Hurtado, D.A. Alvarez, J.A. Paredes,
\newblock {Interval reliability analysis under the specification of statistical information on the input variables},
\newblock {Structural Safety}, 65 (2017) 35--48.
https://doi.org/10.1016/j.strusafe.2016.12.005.

\bibitem{2018:jiang}
C. Jiang, J. Zheng, X. Han,
\newblock {Probability-interval hybrid uncertainty analysis for structures with both aleatory and epistemic uncertainties: a review},
\newblock {Structural and Multidisciplinary Optimization}, 57 (2018) 2485--2502.
https://doi.org/10.1007/s00158-017-1864-4.

\bibitem{1997:kallenberg}
O. Kallenberg,
\newblock {Foundations of Modern Probability},
\newblock Springer, New York, 1997.

\bibitem{2021:MOKarakasevic}
J. Karaka\v{s}evi\'{c}, M. Oberguggenberger,
\newblock Random set solutions to partial differential equations,
\newblock in: A. Sofi, G. Muscolino, R.L. Muhanna (Eds.) {REC 2021: Proceedings of the 9th International Workshop on Reliable Engineering Computing}, 2021,
pp. 279--294. http://ww2new.unime.it/REC2021/proceedings/REC2021\_Proceedings.pdf.

\bibitem{1976:kato}
T. Kato,
\newblock {Perturbation Theory for Linear Operators},
\newblock Springer-Verlag, Berlin-New York, 1976.

\bibitem{2022:kitahara}
M. Kitahara, S. Bi, M. Broggi, M. Beer,
\newblock {Nonparametric {B}ayesian stochastic model updating with hybrid uncertainties},
\newblock {Mechanical Systems and Signal Processing}, 163 (2022) 108195.
https://doi.org/10.1016/j.ymssp.2021.108195.

\bibitem{2020:liu}
H.B. Liu, C. Jiang, Z. Xiao,
\newblock {Efficient uncertainty propagation for parameterized p-box using sparse-decomposition-based polynomial chaos expansion},
\newblock {Mechanical Systems and Signal Processing}, 138 (2020) 106589.
https://doi.org/10.1016/j.ymssp.2019.106589.

\bibitem{2018:lue}
H. L\"{u}, W.-B. Shangguan, D. Yu,
\newblock {A new hybrid uncertainty analysis method and its application to squeal analysis with random and interval variables},
\newblock {Probabilistic Engineering Mechanics}, 51 (2018) 1--10.
https://doi.org/10.1016/j.probengmech.2017.11.001.

\bibitem{2008:matthies}
H.G. Matthies,
\newblock {Stochastic finite elements: computational approaches to stochastic partial differential equations},
\newblock {ZAMM. Zeitschrift f\"{u}r Angewandte Mathematik und Mechanik}, 88 (2008) 849--873.
https://doi.org/10.1002/zamm.200800095.

\bibitem{2016:mishra}
S. Mishra, Ch. Schwab, J. \v{S}ukys,
\newblock {Multi-level {M}onte {C}arlo finite volume methods for uncertainty quantification of acoustic wave propagation in random heterogeneous layered medium},
\newblock {Journal of Computational Physics}, 312 (2016) 192--217.
https://doi.org/10.1016/j.jcp.2016.02.014.

\bibitem{2005:molchanov}
I. Molchanov,
\newblock {Theory of Random Sets},
\newblock Springer-Verlag, Berlin, 2005.

\bibitem{2016:muscolino}
G. Muscolino, R. Santoro, A. Sofi,
\newblock {Reliability analysis of structures with interval uncertainties under stationary stochastic excitations},
\newblock {Computer Methods in Applied Mechanics and Engineering}, 300 (2016) 47--69.
http://dx.doi.org/10.1016/j.cma.2015.10.023.

\bibitem{2013:muscolino}
G. Muscolino, A. Sofi, M. Zingales,
\newblock {One-dimensional heterogeneous solids with uncertain elastic modulus
in presence of long-range interactions: Interval versus stochastic analysis},
\newblock {Computers \& Structures}, 122 (2013) 217--229.
http://dx.doi.org/10.1016/j.compstruc.2013.03.005.

\bibitem{1967:necas}
J. Ne\v{c}as,
\newblock {Les M\'{e}thodes Directes en Th\'{e}orie des \'{E}quations Elliptiques},
\newblock Masson et Cie, Paris; Academia, Prague, 1967.

\bibitem{2020:nedeljkovic}
J. Nedeljkovi\'{c},
\newblock {Partial differential equations with random set and random field coefficients},
\newblock PhD thesis, Universit\"{a}t Innsbruck, Austria, 2020.
https://diglib.uibk.ac.at/urn:nbn:at:at-ubi:1-68791.

\bibitem{2006:nguyen}
H.T. Nguyen,
\newblock {An Introduction to Random Sets},
\newblock Chapman \& Hall, Boca Raton, 2006.

\bibitem{MOBook}
M. Oberguggenberger,
\newblock {Multiplication of Distributions and Applications to Partial Differential Equations},
\newblock Longman, Harlow UK, 1992.

\bibitem{2013:oberguggenberger}
M. Oberguggenberger,
\newblock Combined methods in nondeterministic mechanics,
\newblock in: I. Elishakoff, C. Soize (Eds.), {Nondeterministic Mechanics},
Springer, Wien, 2013, pp. 263--356.

\bibitem{2014:oberguggenberger}
M. Oberguggenberger,
\newblock Stochastic response surfaces, interval analysis, and the reliability of structures,
\newblock in:  M. Beer, S.-K. Au, J.W. Hall (Eds.), {Vulnerability, Uncertainty, and Risk: Quantification, Mitigation and Management},
American Society of Civil Engineers, Reston VA, 2014, pp. 864--875.
http://dx.doi.org/10.1061/9780784413609.

\bibitem{2015:oberguggenberger}
M. Oberguggenberger,
\newblock {Analysis and computation with hybrid random set stochastic models},
\newblock {Structural Safety}, 52 (2015) 233--243.
http://dx.doi.org/10.1016/j.strusafe.2014.05.008.

\bibitem{2008:MOFellin}
M. Oberguggenberger, W. Fellin,
\newblock {Reliability bounds through random sets: Non-parametric methods and geotechnical applications},
\newblock {Computers \& Structures}, 86 (2008) 1093--1101.
https://doi.org/10.1016/j.compstruc.2007.05.040.

\bibitem{2020:MOSchwarz}
M. Oberguggenberger, M. Schwarz,
\newblock On the measurability of stochastic {F}ourier integral operators,
\newblock in: P. Boggiatto, M. Cappiello, E. Cordero, S. Coriasco, G. Garello,
A. Oliaro, J. Seiler (Eds.), {Advances in Microlocal and Time-Frequency Analysis},
Birkh\"{a}user/Springer, Cham, 2020, pp. 383--401.  .

\bibitem{2019:wurzer}
M. Oberguggenberger, L. Wurzer,
\newblock {Random set solutions to stochastic wave equations},
\newblock {Proceedings of Machine Learning Research}, 103 (2019) 314--323.
https://proceedings.mlr.press/v103/oberguggenberger19a.html.

\bibitem{2020:sadeghi}
J. Sadeghi, M. de Angelis, E. Patelli,
\newblock {Robust propagation of probability boxes by interval predictor models},
\newblock {Structural Safety}, 82 (2020) 101889.
https://doi.org/10.1016/j.strusafe.2019.101889.

\bibitem{2010:schmelzer}
B. Schmelzer,
\newblock On solutions to stochastic differential equations with parameters modeled by random sets,
\newblock {International Journal of Approximate Reasoning}, 51 (2010) 1159--1171.
https://doi.org/10.1016/j.ijar.2010.08.006.

\bibitem{2013:schmelzer}
B. Schmelzer,
\newblock Set-valued assessments of solutions to stochastic differential equations with random set parameters,
\newblock {Journal of Mathematical Analysis and Applications}, 400 (2013) 425--438.
https://doi.org/10.1016/j.jmaa.2012.11.019.

\bibitem{2015a:schmelzer}
B. Schmelzer,
\newblock Sklar's theorem for minitive belief functions,
\newblock {International Journal of Approximate Reasoning}, 63 (2015) 48--61.
https://doi.org/10.1016/j.ijar.2015.05.010.

\bibitem{2015:schmelzer}
B. Schmelzer,
\newblock Joint distributions of random sets and their relation to copulas,
\newblock {International Journal of Approximate Reasoning}, 65 (2015) 59--69.
https://doi.org/10.1016/j.ijar.2015.01.007.

\bibitem{2019:schmelzer}
B. Schmelzer,
\newblock Multivariate capacity functionals vs. capacity functionals on product spaces,
\newblock {Fuzzy Sets and Systems}, 364 (2019) 1--35.
https://doi.org/10.1016/j.fss.2018.07.005.

\bibitem{2010:schmelzer:adam}
B. Schmelzer, M. Oberguggenberger,  C.~Adam,
\newblock Efficiency of tuned mass dampers with uncertain parameters on the
  performance of structures under stochastic excitation,
\newblock {Journal of Risk and Reliability}, 224 (2010) 297--308.
https://doi.org/10.1243/1748006XJRR310.

\bibitem{2017:schoebi}
R. Sch\"{o}bi, B. Sudret,
\newblock {Uncertainty propagation of p-boxes using sparse polynomial chaos expansions},
\newblock {Journal of Computational Physics}, 339 (2017) 307--327.
https://doi.org/10.1016/j.jcp.2017.03.021.

\bibitem{2015:sofi}
A. Sofi,
\newblock {Structural response variability under spatially dependent uncertainty:
Stochastic versus interval model},
\newblock {Probabilistic Engineering Mechanics}, 42 (2015) 78--86.
http://dx.doi.org/10.1016/j.probengmech.2015.09.001

\bibitem{2022:Sofi}
A. Sofi,  F. Giunta, G. Muscolino,
\newblock {Reliability analysis of randomly excited FE modelled structures with interval mass and stiffness via sensitivity analysis},
\newblock {Mechanical Systems and Signal Processing}, 163 (2022) 107990.
https://doi.org/10.1016/j.ymssp.2021.107990.

\bibitem{2020:Sofi}
A. Sofi, G. Muscolino, F. Giunta,
\newblock {Propagation of uncertain structural properties described by imprecise Probability Density Functions via response surface method},
\newblock {Probabilistic Engineering Mechanics}, 60 (2020) 103020.
https://doi.org/10.1016/j.probengmech.2020.103020.

\bibitem{2018:sofi:romeo}
A. Sofi, E. Romeo,
\newblock {A unified response surface framework for the interval and stochastic finite element analysis of structures with uncertain parameters},
\newblock {Probabilistic Engineering Mechanics}, 54 (2018) 25--36.
http://dx.doi.org/10.1016/j.probengmech.2017.06.004.

\bibitem{2020:song}
J. Song, M. Valdebenito, P. Wei, M. Beer, Z. Lu,
\newblock {Non-intrusive imprecise stochastic simulation by line sampling},
\newblock {Structural Safety}, 84 (2020) 101936.
https://doi.org/10.1016/j.strusafe.2020.101936.

\bibitem{1987:stroock}
D.W. Stroock,
\newblock {Lectures on Stochastic Analysis: Diffusion Theory},
\newblock Cambridge University Press, Cambridge, 1987.

\bibitem{2004:tonon}
F. Tonon,
\newblock {Using random set theory to propagate epistemic uncertainty through a mechanical system},
\newblock {Reliability Engineering \& System Safety}, 85 (2004) 169--181.
https://doi.org/10.1016/j.ress.2004.03.010.

\bibitem{2010:tonon}
F. Tonon, C.L. Pettit,
\newblock {Toward a definition and understanding of correlation for variables constrained by random relations},
\newblock {International Journal of General Systems}, 39 (2010) 577--604.
https://doi.org/10.1080/03081070903541273.

\bibitem{1975:treves}
F. Tr\`eves,
\newblock {Basic Linear Partial Differential Equations},
\newblock Academic Press, New York, 1975.

\bibitem{2017:troffaes}
M.C.M. Troffaes,
\newblock {A note on imprecise {M}onte {C}arlo over credal sets via importance sampling},
\newblock {Proceedings of Machine Learning Research}, 62 (2017) 325--332.
https://proceedings.mlr.press/v62/troffaes17a.html.

\bibitem{2018:troffaes}
M.C.M. Troffaes,
\newblock {Imprecise Monte Carlo simulation and iterative importance sampling for the estimation of lower previsions},
\newblock {Mechanical Systems and Signal Processing}, 101 (2018) 31--48.
https://doi.org/10.1016/j.ijar.2018.06.009.

\bibitem{2018:troffaesfetzMO}
M.C.M. Troffaes, T. Fetz, M. Oberguggenberger,
\newblock Iterative importance sampling for estimating expectation bounds under partial probability specifations,
\newblock in: M. De Angelis (Ed.), {REC 2018: Proceedings of the 8th International Workshop on Reliable Engineering Computing},
Liverpool University Press, Liverpool, 2002, pp. 147--154.
http://www.rec2018.uk/papers/proceedings/proceedings.pdf.

\bibitem{2013:verhaeghe}
W. Verhaeghe, W. Desmet, D. Vandepitte, D. Moens,
\newblock {Interval fields to represent uncertainty on the output side of a static FE analysis},
\newblock {Computer Methods in Applied Mechanics and Engineering}, 260 (2013) 50--62.
https://doi.org/10.1016/j.cma.2013.03.021.

\bibitem{2019a:wei}
P. Wei, J. Song, S. Bi, M. Broggi, M. Beer, Z. Lu, Z. Yue,
\newblock {Non-intrusive stochastic analysis with parameterized imprecise probability models: {I}. {P}erformance estimation},
\newblock {Mechanical Systems and Signal Processing}, 124 (2019) 349--368.
https://doi.org/10.1016/j.ymssp.2019.01.058.

\bibitem{2019b:wei}
P. Wei, J. Song, S. Bi, M. Broggi, M. Beer, Z. Lu, Z. Yue,
\newblock {Non-intrusive stochastic analysis with parameterized imprecise probability models: {II}. {R}eliability and rare events analysis},
\newblock {Mechanical Systems and Signal Processing}, 126 (2019) 227--247.
https://doi.org/10.1016/j.ymssp.2019.02.015.

\bibitem{2020:zhang}
H. Zhang,
\newblock {Interval importance sampling method for finite element-based structural reliability assessment under parameter uncertainties},
\newblock {Structural Safety}, 38 (2012) 1--10.
https://doi.org/10.1016/j.strusafe.2012.01.003.

\bibitem{2018:zhang}
J. Zhang, M.D. Shields,
\newblock {On the quantification and efficient propagation of imprecise probabilities resulting from small datasets},
\newblock {Mechanical Systems and Signal Processing}, 98 (2018) 465--483.
http://dx.doi.org/10.1016/j.ymssp.2017.04.042.

\end{thebibliography}
\end{document}